\documentclass[12pt]{amsart}
\usepackage{amssymb}  
\usepackage{latexsym} 
\usepackage[all]{xy}
\usepackage{comment}
\usepackage{graphicx}
\usepackage{url}

\newcommand{\C}{{\mathbb C}}
\newcommand{\la}{{\lambda}}

\newcommand{\CC}{{\mathcal C}}
\newcommand{\OO}{{\mathcal O}}
\newcommand{\NN}{{\mathcal N}}
\newcommand{\Char}{\operatorname{char}}

\newcommand{\Diag}{\operatorname{Diag}}
\newcommand{\adj}{\operatorname{adj}}

\newcommand{\dv}{{\,\mid\,}}
\newcommand{\eps}{{\varepsilon}}

\newcommand{\G}{{\mathcal G}}

\newcommand{\GL}{\operatorname{GL}}
\newcommand{\Jac}{\operatorname{Jac}}

\newcommand{\kbar}{\overline{k}}

\newcommand{\Magma}{{\sf MAGMA }}

\newcommand{\nr}{{\text{\rm nr}}}

\newcommand{\OK}{{\mathcal{O}_K}}
\newcommand{\pf}{\operatorname{pf}}

\newcommand{\rd}{\operatorname{rd}}

\newcommand{\PP}{{\mathbb P}}
\newcommand{\Aff}{{\mathbb A}}
\newcommand{\isom}{{\, \cong \,}}
\newcommand{\Q}{{\mathbb Q}}
\newcommand{\rdet}{\operatorname{rd}}
\newcommand{\R}{{\mathbb R}}

\newcommand{\rank}{\operatorname{rank}}
\newcommand{\ra}{{\longrightarrow}}

\newcommand{\x}{{\bf x}}

\newcommand{\Z}{{\mathbb Z}}

\newenvironment{Proof}{\par\noindent{\sc Proof:}}%
                      {\hspace*{\fill}\nobreak$\Box$\par\medskip}
                       {\hspace*{\fill}\nobreak$\Box$\par\medskip}

\newtheorem{Proposition}{Proposition}[section]
\newtheorem{Theorem}[Proposition]{Theorem}
\newtheorem{Lemma}[Proposition]{Lemma}
\newtheorem{Corollary}[Proposition]{Corollary}

\theoremstyle{definition}
\newtheorem{Definition}[Proposition]{Definition}
\newtheorem{Remark}[Proposition]{Remark}
\newtheorem{Example}[Proposition]{Example}
\newtheorem{Algorithm}[Proposition]{Algorithm}

\numberwithin{equation}{section}

\newcommand{\Sha}{\mbox{\wncyr Sh}}
\newfont{\wncyr}{wncyr10 at 12pt}

\begin{document}

\title[Height bounds for coverings of elliptic curves]%
{Local solubility and height bounds for coverings of elliptic curves}


\author{T.A.~Fisher}
\address{University of Cambridge,
         DPMMS, Centre for Mathematical Sciences,
         Wilberforce Road, Cambridge CB3 0WB, UK}
\email{T.A.Fisher@dpmms.cam.ac.uk}

\author{G.F.~Sills}
\address{University of Cambridge,
         DPMMS, Centre for Mathematical Sciences,
         Wilberforce Road, Cambridge CB3 0WB, UK}
\email{gs300@cantab.net}

\date{23rd March 2011} 

\begin{abstract}
We study genus one curves that arise as $2$-, $3$- and $4$-coverings
of elliptic curves. We describe efficient algorithms for testing local
solubility and modify the classical 
formulae for the covering maps so that they work in all 
characteristics. These ingredients are then combined to give
explicit bounds relating the height of a rational point on one of the 
covering curves to the height of its image on the elliptic curve. 
We use our results to improve the existing methods
for searching for rational points on elliptic curves.
\end{abstract}

\maketitle

\renewcommand{\baselinestretch}{1.1}
\renewcommand{\arraystretch}{1.3}
\renewcommand{\theenumi}{\roman{enumi}}

\section{Introduction}

Let $E$ be an elliptic curve over a number field $K$. 
An $n$-covering of $E$ is a
smooth curve of genus one $\CC$ together with a morphism $\pi : \CC \to E$,
with $\CC$ and $\pi$ both defined over $K$, such that the diagram 
\[ \xymatrix{ \CC \ar[d]_{\psi} \ar[dr]^\pi \\ E \ar[r]_{[n]} & E } \]
commutes for some isomorphism $\psi : \CC \isom E$ defined over 
$\overline{K}$. An $n$-descent calculation computes equations for
the everywhere locally soluble $n$-coverings of $E$, i.e. the $n$-coverings 
$\CC$ with $\CC(K_v) \not= \emptyset$ for all places $v$ of $K$. 
Finding  rational points on these $n$-coverings can assist in computing
generators for the Mordell-Weil group $E(K)$. Indeed 
if $\CC(K)$ is non-empty then $\pi (\CC(K))$ is a 
coset of $nE(K)$ in $E(K)$. 

Suppose that $\CC$ is everywhere locally soluble. 
By \cite[Proof of Theorem 1.3]{CasselsIV} there exists 
a $K$-rational divisor $D$ on $\CC$ with 
$D \sim \psi^*(n.\OO)$, where $\OO$ is the identity on $E$.
The complete linear system $|D|$ defines a morphism $\CC \to \PP^{n-1}$.
If $n=2$ then $\CC \to \PP^1$ is a double cover ramified at $4$ points.
If $n \ge 3$ then $\CC \subset \PP^{n-1}$ is a genus one normal curve of
degree $n$. The map $\pi : \CC \to E$ may be recovered as 
$P \mapsto [nP-D] \in {\operatorname{Pic}}^0(\CC) = E$ 
where $D$ is now the hyperplane
section on $\CC$. In the cases $n=2,3,4$ equations for $\CC$ 
take the form of a binary quartic, ternary cubic or quadric intersection.
The Jacobian elliptic curve $E$ and covering map $\pi$ are then 
given by formulae from classical invariant theory 
as surveyed in \cite{AKM3P}.

It is expected that points on $\CC(K)$ will be smaller (and hence
easier to find) than their images in $E(K)$. This statement is
made precise using the theory of heights. 
Let $h$ be the logarithmic height on $\CC$
relative to the hyperplane section $D$, and $h_E$ the $x$-coordinate
logarithmic height on $E$. Then as pointed out in \cite{StollParis}
there exist constants $B_1$ and $B_2$ such that 
\begin{equation}
\label{eqn1}
 B_1 \le h(P) - \frac{1}{2n} h_E(\pi P) \le B_2 
\end{equation}
for all $P \in \CC(K)$.
To prove this one first notes that since $n^2. \OO \sim [n]^* \OO$ 
we have $2n D \sim \pi^* (2.\OO)$. The existence of bounds $B_1$ and $B_2$
then follows by standard results about heights; 
see for example \cite[Theorem B.3.2]{HS}.

We restrict to $n=2,3$ or $4$. In these cases $n$-descent has been 
implemented in the computer algebra system Magma~\cite{Magma}
at least over $K= \Q$. The algorithms for $3$-descent are described
in \cite{SchaeferStoll}, \cite{ndescent} and those for $4$-descent
in \cite{MSS}, \cite{WomackThesis}.
In Sections~\ref{sec:locsol},~\ref{sec:covmaps} 
and~\ref{sec:htbds} we
\begin{itemize}
\item describe algorithms for testing whether $\CC(K_v) \not = \emptyset$,
\item modify the formulae for the covering map $\pi : \CC \to E$ 
so that they work in all characteristics, and
\item compute explicit bounds $B_1$ and $B_2$ in (\ref{eqn1}).
\end{itemize}


Recent work on higher descents and on computing the Cassels-Tate pairing 
(see \cite{CreutzThesis}, \cite{Donnelly}, \cite{cbrank}, 
\cite{StammingerThesis}) relies on being able to efficiently compute 
local points. This prompted us to improve the local
solubility tests currently implemented in Magma. The material
in Section~\ref{sec:locsol} should however contain few surprises for 
experts. The main reason for including it here is as a preliminary to 
our work on height bounds. The latter is also the motivation for the
formulae in Section~\ref{sec:covmaps},
although these too may be of independent interest.

It is possible to compute bounds $B_1$ and $B_2$ in (\ref{eqn1})
using elimination theory. However this method gives rather poor results.
Instead we compute our bounds as sums of local contributions. 
This generalises work of Siksek \cite{Siksek}
who considered the case where $\pi$ is multiplication-by-$2$ on $E$.
As he observes it is worth putting some effort into obtaining good
bounds, as this can significantly reduce the size of the region
we end up searching. We give some examples in Section~\ref{sec:examples}.

The bounds $B_1$ and $B_2$ depend on our choice of equations for 
$\CC$ and $E$. Let us take $K = \Q$. For $E$ we take the global
minimal Weierstrass equation
\begin{equation}
\label{minw}
 y^2 + a_1 x y + a_3 y = x^3 + a_2 x^2 + a_4 x + a_6 
\end{equation}
with $a_1,a_3 \in \{0,1\}$ and $a_2 \in \{0,\pm 1\}$. For $\CC$ we take
an equation that is minimised and reduced as described in \cite{CFS}.
Roughly speaking one expects that minimising improves the bounds
at the finite places, and reducing improves the bounds at the
infinite places. However there can be more than one choice of
minimisation. We find that the bounds can vary significantly between 
these choices. In Section~\ref{sec:examples} we include an example
where these ideas allow us to improve the search for rational points
on $\CC$ (and hence on $E$).

\subsection{Genus one models}
\label{sec:g1m}
The following notation is recalled from \cite{CFS}, \cite{g1inv}.
We call the equations defining an $n$-covering (where $n=2,3$ or $4$)
a {\em genus one model}. 
More precisely we make the following definition.

\begin{Definition}
\label{def:g1m}
Let $R$ be any ring. 
\begin{enumerate}
\item A {\em genus one model of degree $2$} over $R$ is a generalised binary 
quartic
\[ y^2 + P(x_1,x_2) y = Q(x_1,x_2), \]
sometimes abbreviated $(P,Q)$, 
where $P$ and $Q$ are homogeneous forms of degree 2 and 4 with 
coefficients in $R$. A transformation of genus one models is given by 
$y \leftarrow \mu^{-1} y + r_0 x_1^2 + r_1 x_1 x_2 + r_2 x_2^2$ for 
some $\mu \in R^\times$ and $r=(r_0,r_1,r_2) \in R^3$, followed by
$x_j \leftarrow \sum n_{ij} x_i$ for some $N = (n_{ij}) \in \GL_2(R)$.
We write $\G_2(R)$ for the group of all such transformations 
$g= [\mu,r,N]$ and define $\det g = \mu \det N$.
\item A {\em genus one model of degree $3$} over $R$ is a ternary cubic 
$U \in R[x_1,x_2,x_3]$. A transformation of genus one models is given 
by multiplying the cubic through by $\mu \in R^\times$, followed by
$x_j \leftarrow \sum n_{ij} x_i$ for some $N = (n_{ij}) \in \GL_3(R)$. 
We write $\G_3(R)$ for the group of all such transformations 
$g= [\mu,N]$ and define $\det g = \mu \det N$.
\item A {\em genus one model of degree $4$} over $R$ is a quadric intersection,
i.e. a pair of homogeneous polynomials $Q_1,Q_2 \in R[x_1, \ldots, x_4]$  
of degree $2$. A transformation of quadric intersections is given by
$Q_i \leftarrow \sum m_{ij} Q_j$ for some $M = (m_{ij}) \in \GL_2(R)$ 
and $x_j \leftarrow \sum n_{ij} x_i$ for some $N = (n_{ij}) \in \GL_4(R)$.
We write $\G_4(R)$ for the group of all such transformations 
$g= [M,N]$ and define $\det g = \det M \det N$.
\end{enumerate}
\end{Definition} 
 
We say that genus one models are {\em $R$-equivalent} if they are 
in the same orbit for the action of $\G_n(R)$. Notice that by our 
conventions the action of $\G_n(R)$ on the space of genus one models 
is a left action.

An invariant of weight $k$ is a polynomial $F$ in the coefficients of
a genus one model such that $F \circ g = (\det g)^k F$ for all 
$g \in \G_n$. Let $c_4, c_6$ and $\Delta=(c_4^3-c_6^2)/1728$ 
be the classical invariants of weights $4$, $6$ and $12$. 
We fix the scaling of these invariants as described in 
\cite{CFS}, \cite{g1inv}, i.e so that the models $y^2 + x_1 x_2 y = 0$,
$x_1 x_2 x_3 = 0$ and  $x_1 x_2 = x_3 x_4=0$ have invariants 
$c_4 = 1$ and $c_6 = -1$. For example the binary quartic
$y^2 = a x^4 + b x^3 z + c x^2 z^2 + d x z^3 + e z^4$ 
has invariants
\begin{align*}
 c_4  & = 2^4(12 a e - 3 b d + c^2) \\
 c_6  & = 2^5(72 a c e - 27 a d^2 - 27 b^2 e + 9 b c d - 2 c^3). 
\end{align*}

A genus one model $\Phi$ over a field $K$ is {\em non-singular} 
if the variety $\CC_{\Phi}$ it 
defines is a smooth curve of genus one, and {\em $K$-soluble} if 
$\CC_{\Phi}(K) \not= \emptyset$. It is shown in \cite{g1inv}
that $\Phi$ is non-singular if and only if $\Delta(\Phi) \not = 0$. 
Moreover if $\Char(K) \not=2,3$ then
(by an observation originally due to Weil in the cases $n=2,3$)
the Jacobian elliptic curve $E = \Jac(\CC_{\Phi})$ has Weierstrass equation
\begin{equation}
\label{eqn:jac}
  y^2 = x^3 - 27 c_4(\Phi)x - 54 c_6 (\Phi). 
\end{equation}
Functions for computing with genus one models, their transformations
and invariants have been contributed to Magma \cite{Magma} 
by the first author.

\section{Testing for local solubility}
\label{sec:locsol}

Let $K$ be a finite extension of $\Q_p$ with ring of integers
$\OK$, maximal ideal $\pi \OK$, residue field $k$ and normalised discrete 
valuation $v : K^\times \to \Z$. Reduction mod $\pi$ will be denoted 
$x \mapsto \widetilde{x}$. If $f$ is a polynomial with coefficients
in $K$ then we write $v(f)$ for the minimum valuation of a coefficient.

Let $\Phi$ be a non-singular genus one model over $K$ of degree $n \in
\{2,3,4\}$. In this section we give algorithms for deciding whether
$\Phi$ is $K$-soluble.  Our algorithm in the case $n=2$ is essentially
the same as that in \cite{BSD1}, \cite{Bruin}, \cite{Cremonabook}, 
\cite{MSS} and is
included only for completeness. The cases $n=3,4$ can also be handled
by the general method for complete intersections described in
\cite{Bruin}.  However this general method involves looping over all
$k$-points on the reduction, and is therefore inefficient when 
$k$ is large. We overcome this problem by making use of the 
geometry of singular genus one models. We have contributed our algorithms
(over $K= \Q_p$) to Magma \cite{Magma}, and from the 
next release (Version 2.17) they 
will be called by default when equations of the relevant form are 
passed to {\tt IsLocallySoluble}.

The basic algorithms are listed in Section~\ref{sec:algs}. They 
depend on methods for deciding whether there are any smooth $k$-points 
on the reduction (see Section~\ref{sec:smooth}) and for 
finding all non-regular 
$k$-points (see Section~\ref{sec:nonreg}). It is clear
by Hensel's lemma that when an answer is returned 
then that answer is correct. If the algorithms failed to terminate then
from the resulting infinite sequence of transformations we could 
construct a singular point on the original curve. Thus our
assumption that $\Phi$ is non-singular ensures that the algorithms
terminate. We omit the details since we give an alternative
proof in Section~\ref{sec:bounds}.

In practice we first replace $\Phi$ by a {\em minimal model}, i.e. 
a $K$-equivalent model over $\OK$ with $v(\Delta(\Phi))$ minimal. 
Algorithms for doing this are described in \cite{CFS}.
Let $E = \Jac(\CC_{\Phi})$ be the Jacobian elliptic curve 
and $\Delta_E$ its minimal
discriminant. Then $v(\Delta(\Phi)) = v (\Delta_E) + 12 \ell$ where
$\ell$ is a non-negative integer called the {\em level} of $\Phi$. 
Notice that applying a transformation $g \in \G_n(K)$ changes the
level by $v(\det g)$.
In \cite{CFS} it is shown that the minimal level is $0$ if and 
only if $\CC_{\Phi}(K^\nr) \not=\emptyset$ where $K^\nr$ is the 
maximal unramified extension of $K$.
Therefore our local solubility tests are only needed for models of 
level $0$. This extra hypothesis will be useful in Section~\ref{sec:nonreg}.

We mention as an aside that if the Tamagawa number $c(E)$ is coprime
to $n$ then a further simplification is possible. Indeed by the following
lemma we have $\CC_{\Phi}(K) \not=\emptyset$ if and only if 
$\CC_{\Phi}(K^\nr) \not=\emptyset$, and so the algorithms in \cite{CFS}
already give a test for local solubility. 

\begin{Lemma}
\label{lem:kerres}
The restriction map $H^1(K,E) \to H^1(K^\nr,E)$ has kernel of order $c(E)$.
\end{Lemma}
\begin{Proof}
We recall the argument from the proof of \cite[Theorem 3.1]{AS}.
By \cite[Proposition 3.8]{MilneADT} and 
the inflation-restriction exact sequence 
the kernel is isomorphic to $H^1(k,\Phi_E)$ where
$\Phi_E$ is the component group of the N\'eron model of $E$.
Since $\Phi_E$ is finite and $c(E) = \# \Phi_E(k)$ the result
follows by the exact sequence
\[ 0 \ra H^0(k,\Phi_E) \ra \Phi_E \stackrel{1-{\rm Frob}}{\ra}
\Phi_E \ra H^1(k,\Phi_E) \ra 0. \vspace{-0.3cm} \]
\end{Proof}

\enlargethispage{0.2cm}

\subsection{Algorithms}
\label{sec:algs}
Let $\Phi$ be a non-singular genus one model over $K$ 
of degree $n \in \{2,3,4\}$. Our algorithms for deciding whether
$\CC_{\Phi}(K) \not= \emptyset$  start by making two simplifications.
First by clearing denominators 
we may assume that $\Phi$ is defined over $\OK$.
Then by calling the algorithm $n$ times  (with the co-ordinates 
permuted) it suffices to look for points on a standard affine piece
with co-ordinates in $\OK$.
We remark that if $\Char(k) \not= 2$ then the
first algorithm simplifies in the obvious way by completing the square.

\begin{Algorithm}
\label{alg:n=2}
{\tt IsLocallySoluble(h,g) } \\
{\tt INPUT:} Polynomials $h(x),g(x) \in \OK[x]$ with $\deg(h) \le 2$ 
and $\deg(g) \le 4$. \\
{\tt OUTPUT:  TRUE/FALSE} (solubility of $y^2 + h(x) y = g(x)$ 
for $x,y \in \OK$) 
\begin{enumerate}
\item Make a substitution $ y \leftarrow y + r_0 x^2 + r_1 x + r_2$
(with $r_i \in \OK$) so that if possible $v(h) \ge 1$ and $v(g) \ge 1$. 
If now $v(h) \ge 1$ and $v(g) \ge 2$ then replace $h$ and $g$ 
by $\pi^{-1}h$ and $\pi^{-2}g$ and repeat Step (i).
\item Consider the affine curve
\begin{equation*}
\Gamma = \{ y^2 + \widetilde{h}(x) y = \widetilde{g}(x) \} \subset \Aff^2_k . 
\end{equation*}
If there are smooth $k$-points on $\Gamma$ then return {\tt TRUE}.
\item Find all non-regular $k$-points on $\Gamma$.
These are the singular points $(\widetilde{u},\widetilde{v})$ 
on $\Gamma$ with the property that for some (and hence all)
lifts $u,v \in \OK$ of $\widetilde{u},\widetilde{v} \in k$ we have
$v^2 + h(u) v \equiv g(u) \pmod{\pi^2}$. 
\item For each non-regular $k$-point $(\widetilde{u},\widetilde{v})$ 
on $\Gamma$ lift $\widetilde{u} \in k$ to $u \in \OK$ and put 
$h_1(x) = h(u + \pi x)$, $g_1(x) = g(u + \pi x)$.
If {\tt IsLocallySoluble(h1,g1)} then return {\tt TRUE}.
\item Return {\tt FALSE}.
\end{enumerate}
\end{Algorithm}

\begin{Algorithm}
\label{alg:n=3}
{\tt IsLocallySoluble(g) } \\
{\tt INPUT:} A polynomial $g(x,y) \in \OK[x,y]$ of total degree $\le 3$. \\
{\tt OUTPUT:  TRUE/FALSE} (solubility of $g(x,y) = 0$ for $x,y \in \OK$) 
\begin{enumerate}
\item Divide $g$ by $\pi^{v(g)}$ so that now $v(g) = 0$.
\item Consider the affine curve
\begin{equation*}
\Gamma = \{ \widetilde{g}(x,y)=0 \} \subset \Aff^2_k . 
\end{equation*}
If there are smooth $k$-points on $\Gamma$ then return {\tt TRUE}.
\item Find all non-regular $k$-points on $\Gamma$.
These are the singular points $(\widetilde{u},\widetilde{v})$ 
on $\Gamma$ with the property that for some (and hence all)
lifts $u,v \in \OK$ of $\widetilde{u},\widetilde{v} \in k$ we have
$g(u,v) \equiv 0 \pmod{\pi^2}$. 
\item For each non-regular $k$-point $(\widetilde{u},\widetilde{v})$ on $\Gamma$ 
lift $\widetilde{u},\widetilde{v} \in k$ to $u,v \in \OK$ and put 
$g_1(x,y)= g(u + \pi x,v + \pi y)$.
If {\tt IsLocallySoluble(g1)} then return {\tt TRUE}.
\item Return {\tt FALSE}.
\end{enumerate}
\end{Algorithm}

\begin{Algorithm}
\label{alg:n=4}
{\tt IsLocallySoluble(g1,g2) } \\
{\tt INPUT:} Polynomials $g_1,g_2 \in \OK[x,y,z]$ of total degree $\le 2$. \\
{\tt OUTPUT:  TRUE/FALSE} (solubility of $g_1(x,y,z) = g_2(x,y,z) = 0$ for 
$x,y,z \in \OK$) 
\begin{enumerate}
\item Replace $g_1$ and $g_2$ by linear combinations so that
$\widetilde{g}_1$ and $\widetilde{g}_2$ are linearly independent over $k$. 
If $\widetilde{g}_1$ and $\widetilde{g}_2$ have a common linear factor then  
make a change of coordinates so that this factor is $x$. 
Then replace $g_i(x,y,z)$ by $\pi^{-1} g_i(\pi x,y,z)$ for $i=1,2$
and repeat Step (i). 
\item Consider the affine curve
\begin{equation*}
\Gamma = \{ \widetilde{g}_1(x,y,z)=\widetilde{g}_2(x,y,z)=0 \} 
\subset \Aff^3_k . 
\end{equation*}
If there are smooth $k$-points on $\Gamma$ then return {\tt TRUE}.
\item Find all non-regular $k$-points on $\Gamma$.
These are the points $(\widetilde{u},\widetilde{v},\widetilde{w})$
on $\Gamma$ that are singular on $\{ \widetilde{g} = 0 \}$ 
for some $g = \lambda g_1 + \mu g_2$ (where $\lambda,\mu \in \OK$ are not 
both divisible by $\pi$) with the property that for some (and hence all)
lifts $u,v,w \in \OK$ of $\widetilde{u},\widetilde{v},\widetilde{w} \in k$ 
we have $g(u,v,w) \equiv 0 \pmod{\pi^2}$. 
\item For each non-regular $k$-point $(\widetilde{u},\widetilde{v},
\widetilde{w})$ on $\Gamma$ lift $\widetilde{u},\widetilde{v},\widetilde{w} 
\in k$ to $u,v,w \in \OK$ and put 
\[h_i(x,y,z)= g_i(u + \pi x,v + \pi y,w + \pi z)\] for $i=1,2$.
If {\tt IsLocallySoluble(h1,h2)} then return {\tt TRUE}.
\item Return {\tt FALSE}.
\end{enumerate}
\end{Algorithm}

\begin{Remark} 
The algorithms may be adapted to return a certificate in the case $\Phi$ 
is locally soluble. This certificate takes the form of a transformation 
of genus one models $g$ such that $g \Phi$ has smooth 
$k$-points on its reduction. A smooth $k$-point on the reduction is 
easily found (e.g. by intersecting with random hyperplanes). We may then
use Hensel's lemma to compute a local point to any desired precision.
This is the second returned argument of Magma's 
{\tt IsLocallySoluble}.
\end{Remark}

\enlargethispage{0.8cm}

\subsection{Testing for smooth points}
\label{sec:smooth}
We show how to decide whether a genus one
model defined over a finite field $k$ has any smooth $k$-points.
For small $k$ there is no difficulty in looping over all $k$-points 
and testing to see which if any are smooth. For larger $k$ this can be 
rather inefficient.

First we recall the classification of singular genus one
models over an algebraically closed field ${\mathbb K}$.
Notice that we are only interested in models that define a curve.

\begin{Lemma}
\label{cl2}
The $\GL_2({\mathbb K})$-orbits of singular binary quartics have the following 
representatives.
\[ \begin{array}{l@{\qquad}lll}
& \text{ binary quartic } & \text{ geometric description } \\ \hline
A_1 &  y^2 = x^3 z + x^2 z^2  & \text{ a rational nodal curve } \\  
A_2 &  y^2 = x^2 z^2  & \text{ two rational curves } \\  
B_1 &  y^2 = x^3 z  & \text{ a rational cuspidal curve } \\  
B_2 &  y^2 = x^4  & \text{ two rational curves } \\  
D & y^2 = 0 & \text{ a double line } 
\end{array} \]
\end{Lemma}
\begin{Proof}
These cases correspond to the number and multiplicity of the repeated
roots of the binary quartic.
\end{Proof}

\begin{Lemma}
\label{cl3}
Assume $char({\mathbb K}) \not=3$.
Then the $\GL_3({\mathbb K})$-orbits of non-zero singular ternary cubics have the following 
representatives.
\[ \begin{array}{l@{\qquad}lll}
& \text{ ternary cubic } & \text{ geometric description } \\ \hline
A_1 &  xyz - y^3 - z^3  & \text{ a rational nodal cubic } \\  
A_2 &  xyz - y^3  & \text{ a conic and a line } \\  
A_3 &  xyz  & \text{ three lines } \\  
B_1 &  y^2 z - x^3  & \text{ a rational cuspidal cubic } \\  
B_2 &  x^2 y - y^2 z  & \text{ a conic and a line } \\  
B_3 &  x^2 y - x y^2  & \text{ three lines } \\  
C &  x^2 y & \text{ a line and a double line }  \\
D & x^3 & \text{ a triple line } 
\end{array} \]
\end{Lemma}
\begin{Proof}
This is standard. See for example \cite[Section 10.3]{Dolg}.
\end{Proof}

\begin{Lemma}
\label{cl4}
Assume $char({\mathbb K}) \not=2$.
Then the $\GL_2({\mathbb K}) \times \GL_4({\mathbb K})$-orbits of quadric 
intersections $(Q_1,Q_2)$, with
$Q_1$ and $Q_2$ coprime, 
have the following representatives. (The final column relates to
Lemma~\ref{dec4} below.)
\[ \begin{array}{l@{\qquad}llll}
& \text{ quadric intersection } & \text{ geometric description } & 
\text{ \hspace{-2em} Segre symbol} & m \\ \hline
A_1 &  x_1 x_3 - x_2^2 - x_4^2 ,\,\, x_2 x_4 - x_3^2 & \text{ a rational nodal quartic } & [112] & 0 \\  
A_2 &  x_1 x_3 - x_2^2 ,\,\, x_2 x_4 - x_3^2 & \text{ a twisted cubic and a line } & [22] & 0 \\  
A_3 &  x_1 x_4 - x_2^2 - x_3^2 ,\,\, x_2 x_3 & \text{ two conics } & [11(11)] & 1 \\  
A_4 &  x_1 x_3 - x_2^2 ,\,\, x_2 x_4 & \text{ a conic and two lines } & [2(11)] & 1 \\  
A_5 &  x_1 x_3 ,\,\, x_2 x_4 & \text{ four lines } & [(11)(11)] & 2 \\  
B_1 &  x_1 x_4 - x_2^2 ,\,\, x_2 x_4 - x_3^2 & \text{ a rational cuspidal quartic } & [13] & 0 \\  
B_2 &  x_1 x_4 - x_2 x_3 ,\,\, x_2 x_4 - x_3^2 & \text{ a twisted cubic and a line } & [4] & 0 \\  
B_3 &  x_1 x_3 + x_1 x_4 - x_2^2 ,\,\, x_3 x_4 & \text{ two conics } & [1(21)] & 1 \\  
B_4 &  x_1 x_3 - x_2^2 + x_2 x_4 ,\,\, x_3 x_4 & \text{ a conic and two lines } & [(31)] & 1 \\  
B_5 &  x_2 x_3 - x_3 x_4 ,\,\, x_2 x_4 - x_3 x_4 & \text{ four lines } & [111] & 3 \\  
C_1 &  x_2 x_3 - x_3 x_4 ,\,\, x_2 x_4 - x_3 x_4 & \text{ a conic and a double line } & [1\{3\}] & 1 \\  
C_2 &  x_1 x_3 + x_2 x_4 ,\,\, x_1 x_4 & \text{ two lines and a double line } & [(22)] & 1 \\  
C_3 &  x_2 x_3 - x_2 x_4 ,\,\, x_3 x_4 & \text{ two lines and a double line } & [12] & 2 \\  
C_4 &  x_2 x_3 - x_4^2 ,\,\, x_3 x_4 & \text{ a line and a triple line } & [3] & 1 \\  
D_1 &  x_1^2 - x_2 x_3 ,\,\, x_4^2 & \text{ a double conic } & [1(111)] & - \\  
D_2 &  x_1 x_4 + x_2 x_3 ,\,\, x_4^2 & \text{ two double lines } & [(211)] & - \\  
D_3 &  x_2 x_3 ,\,\, x_4^2 & \text{ two double lines } & [1(11)] & - \\  
D_4 &  x_2 x_4 - x_3^2 ,\,\, x_4^2 & \text{ a quadruple line } & [(21)] & - \\  
D_5 &  x_3^2 ,\,\, x_4^2 & \text{ a quadruple line } & [11] & - 
\end{array} \]
\end{Lemma}
\begin{Proof}
The classification (at least over ${\mathbb K} = \C$) is due 
to Segre. See for example \cite{Bromwich}, \cite{DLLP}, \cite{HodgePedoe}.
\end{Proof}

\begin{Remark}
The restrictions on the characteristic of ${\mathbb K}$ in Lemmas~\ref{cl3}
and~\ref{cl4} are necessary. For example if $\Char({\mathbb K}) =3$ 
then the cuspidal cubics $y^2 z = x^3 + \lambda x^2 y$ are 
inequivalent for $\lambda = 0$ and $\lambda \not=0$. 
Likewise if $\Char({\mathbb K}) =2$
then the cuspidal quadric intersections $x_1 x_4 + \lambda x_2 x_3
 - x_2^2 = x_2 x_4 - x_3^2 = 0$  are inequivalent for $\lambda =0$ and
$\lambda \ne 0$. 
\end{Remark}

Figures 1, 2 and 3 illustrate the classifications 
in Lemmas~\ref{cl2}, \ref{cl3} and~\ref{cl4}.

\medskip

\begin{center}
\includegraphics{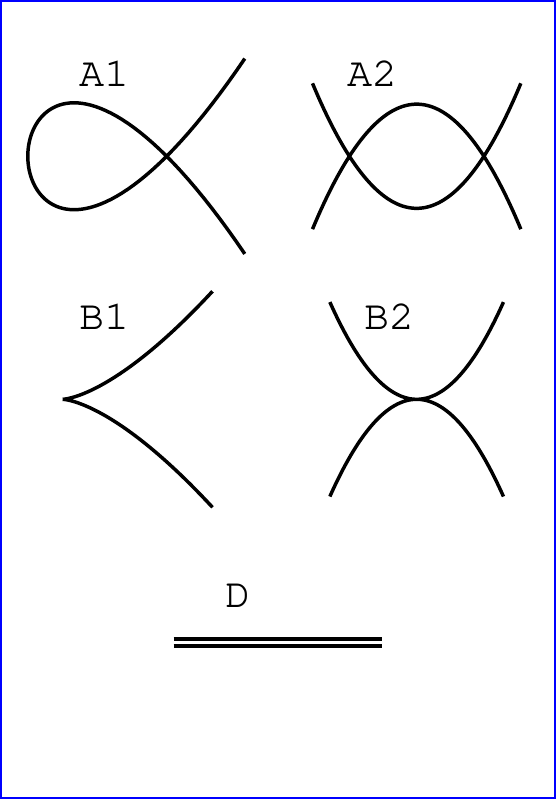}
\includegraphics{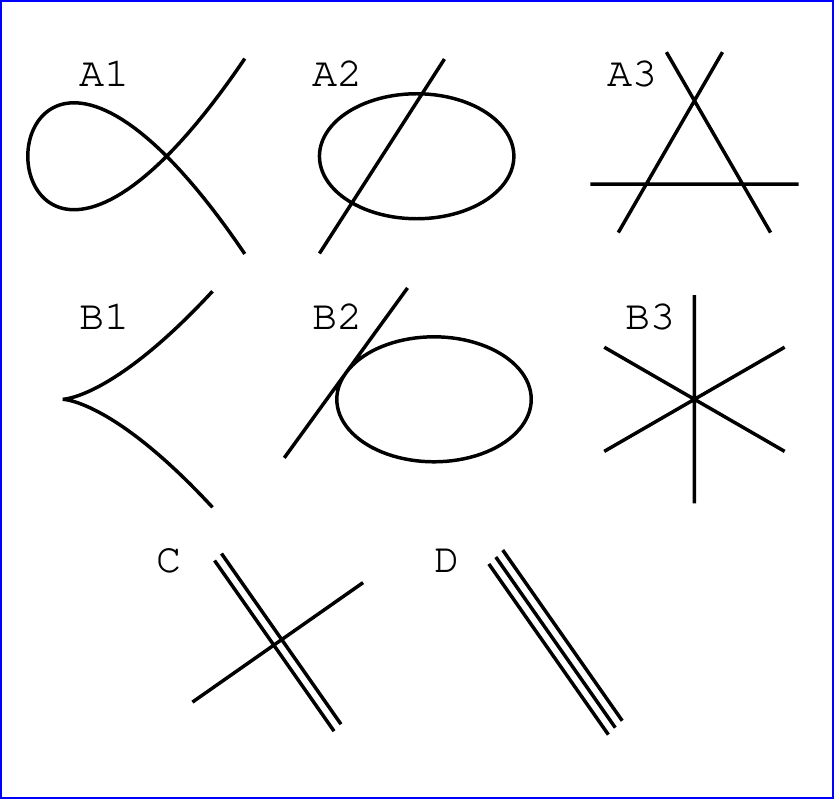}
Figure 1 \hspace{5cm}  Figure 2 \hspace{1cm} ~
\end{center}

\enlargethispage{0.5cm}

\begin{center}
\includegraphics{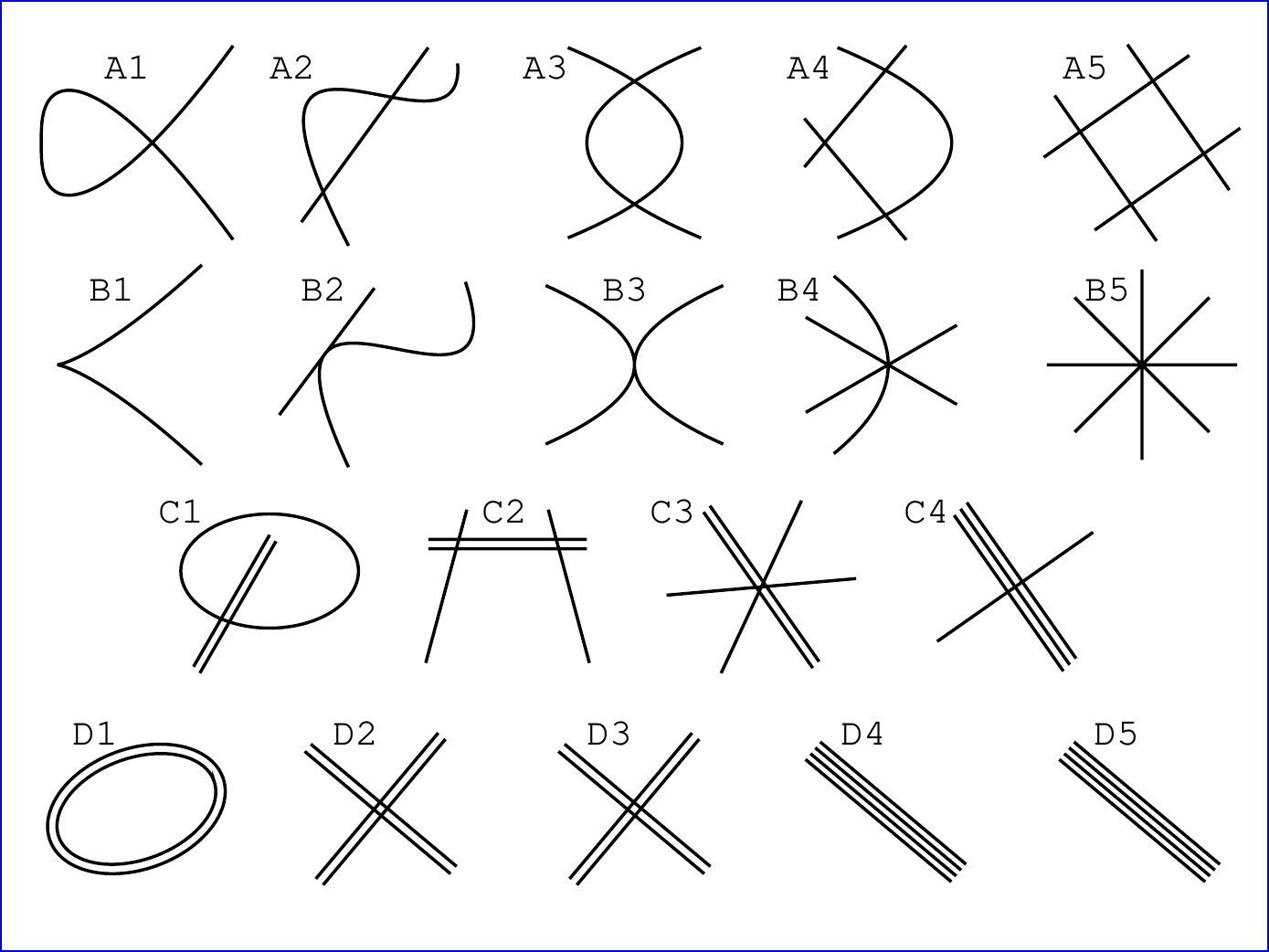}
Figure 3
\end{center}

\medskip

Let $\Phi$ be a genus one model over a finite field $k$.
To decide whether there are any smooth $k$-points on $\CC_\Phi$
we employ the following lemmas.

\begin{Remark}
The algorithms in Section \ref{sec:algs} in fact ask whether
there are any smooth $k$-points on some affine piece $\Gamma$ 
of $\CC_{\Phi}$. It can happen that all the smooth $k$-points lie on
the hyperplane at infinity, either because $k$ is small or because
all relevant components are contained in that hyperplane. 
In terms of our original task of deciding $K$-solubility 
this simply means that we find a point sooner than expected, i.e.
even before we consider an affine piece where it has integral 
co-ordinates.
\end{Remark}

In the case $n=2$ we assume $\Char (k) \not=2$. In particular we may 
complete the square so that our models are given by binary quartics.

\begin{Lemma}
\label{spts2}
Assume $\Char (k) \not= 2$ and let $F \in k[x,z]$ be a binary quartic.
\begin{enumerate}
\item If $F$ is identically zero then $\CC_F$ has no smooth $k$-points.
\item If $F$ is non-zero, but factors as $F(x,z) = \alpha G(x,z)^2$,
then $\CC_F$ has a smooth $k$-point 
if and only if $\alpha \in (k^\times)^2$.
\item In all other cases $\CC_F$ has a smooth $k$-point.
\end{enumerate}
\end{Lemma}

\begin{Proof} 
This is clear by Lemma~\ref{cl2}.
\end{Proof}

We write $\kbar$ for the algebraic closure of $k$.

\begin{Lemma}
\label{cubics}
Let $U \in k[x,y,z]$ be a non-zero ternary cubic.
\begin{enumerate}
\item If $U$ factors over $\kbar$ as a product of linear forms 
then $\CC_U$ has a smooth $k$-point if and only if
one of these linear forms is defined over $k$ and is not a repeated factor.
\item In all other cases $\CC_U$ has a smooth $k$-point.
\end{enumerate}
\end{Lemma}
\begin{Proof} This is clear by Lemma~\ref{cl3}.
\end{Proof}

Now let $\Phi = (Q_1,Q_2)$ be a model of degree $4$.
It is clear that if there is a rank~1 quadric in the pencil 
\begin{equation}
\label{pen}
 \{ \lambda Q_1 + \mu Q_2 \mid (\lambda:\mu) \in \PP^1(\kbar) \} 
\end{equation}
then $\CC_{\Phi}$ has no smooth $k$-points. 

\begin{Lemma}
\label{dec4}
Assume $\Char(k) \not= 2$ and let $\Phi=(Q_1,Q_2)$ be a quadric intersection 
over $k$ with $Q_1$ and $Q_2$ coprime.
Suppose the pencil~(\ref{pen}) over $\kbar$ contains no rank~$1$ quadrics and exactly 
$m$ rank~$2$ quadrics.
\begin{enumerate}
\item If $m=0$ then $\CC_{\Phi}$ has a smooth $k$-point. 
\item If $m=1$ then $\CC_{\Phi}$ has a smooth $k$-point if and only
if the rank $2$ quadric in the pencil factors over $k$. 
\item If $m \ge 2$ then $\CC_{\Phi}$ is (set-theoretically) a union of lines.
\end{enumerate}
\end{Lemma}
\begin{Proof}
This follows from the classification in Lemma~\ref{cl4}.
(The integer $m$ is recorded in the statement of the lemma. 
It is replaced by a dash in cases where there is a rank~$1$ quadric.)
\end{Proof}

It remains to test for smooth $k$-points in the case $\CC_{\Phi}$ 
is a union of lines. Let $A$ and $B$ be 
the 4 by 4 symmetric matrices corresponding to 
$Q_1$ and $Q_2$. Let $M$ be the generic $4$ by $4$ skew-symmetric
matrix. 
The {\em Fano scheme} is the subscheme of $\PP^5$ defined
by the vanishing of the Pfaffian of $M$ and all entries of the
matrices $MAM$ and $MBM$. Identifying $\{ {\operatorname{Pf}}(M) = 0 \} 
\subset \PP^5$ with the Grassmannian of lines in $\PP^3$, 
the points of the Fano scheme
correspond to the lines on the quadric intersection. In particular
the Fano scheme is zero-dimensional.

\begin{Lemma} 
Assume $\Char(k) \not= 2$ and let $\Phi$ be a quadric intersection 
such that $\CC_{\Phi}$ is (set-theoretically) a union of lines.
Then $\CC_\Phi$ has a smooth $k$-point if and only if the Fano scheme 
has a smooth $k$-point.
\end{Lemma}
\begin{Proof}
It suffices to show that a line has multiplicity one if and only if 
it corresponds to a smooth point on the Fano scheme. We checked
this using the classification in Lemma~\ref{cl4}. 
\end{Proof}

\begin{Remark}
Assume $\Char (k) \not=2,3$. 
Then one way to test whether a binary quartic $F$ 
is the square of a polynomial over $\kbar$ is to test whether $F$ and 
its Hessian (which is again a binary quartic) are linearly dependent.
Likewise if $\Phi$ is a genus one model of degree $3$ or $4$ and
$\CC_\Phi$ is a curve then $\CC_\Phi$ is a union of lines if and 
only if $\Phi$ and its Hessian are linearly dependent (as genus one models).
For the definition of the Hessian in the case $n=4$ see \cite{g1hess}.
\end{Remark}

\subsection{Finding the non-regular points}
\label{sec:nonreg}
We keep the notation for local fields introduced at the start of
Section~\ref{sec:locsol}. In particular $K$ is a finite extension 
of $\Q_p$ with ring of integers $\OK$ and residue field $k$.

We show how to find the $k$-rational non-regular points on the 
reduction of a genus one model over $K$.
(See the algorithms of Section~\ref{sec:algs} for the definition 
of a non-regular point.)
If $k$ is small or the singular locus is zero-dimensional 
then there is no difficulty in
looping over all singular points on the reduction and testing
to see which if any are non-regular.
For larger $k$ this can be rather inefficient. Instead we employ 
the following lemmas. 

Recall that by the results in~\cite{CFS} we may
assume that our models have level~$0$ and so in 
particular are minimal. Notice also that, taking into account
the transformations in Step (i) that immediately follow each
recursion, the algorithms in Section~\ref{sec:algs} never increase
the level.

\begin{Lemma}
Assume $\Char(k) \not= 2$ and let $y^2 = F(x,z)$ be a minimal binary
quartic over $K$. Then the non-regular points are
some (but not necessarily all) of the roots 
of $F_1(x,z) \equiv 0 \pmod{\pi}$ where $F_1 = \pi^{-v(F)} F$.
\end{Lemma}
\begin{Proof} Since $F$ is minimal we have $v(F) = 0$ or $1$. 
The rest is clear.
\end{Proof}  

\begin{Lemma}
Let $F(x,y,z)$ be a minimal ternary cubic over $K$.
If the singular locus of the reduction has positive dimension 
then by a change of co-ordinates we may assume that
\[  F(x,y,z) = f_0 x^3 + f_1(y,z) x^2 + \pi f_2(y,z) x + \pi f_3(y,z) \]
where the $f_i$ are binary forms of degree $i$. There are then at most
$3$ non-regular points and these are the roots of 
$x \equiv f_3(y,z) \equiv 0 \pmod{\pi}$.
\end{Lemma}
\begin{Proof} Since $F$ is minimal we have $v(F) = 0$ and $v(f_3)=0$.
The rest is clear.
\end{Proof}  

Assume $\Char(k) \not= 2$ and consider the quadric intersection
$\x^T A \x = \x^T B \x = 0$ where $A = (a_{ij})$ and $B=(b_{ij})$ are
$4$ by $4$ symmetric matrices over $\OK$. Then $(1:0:0:0)$ is a non-regular
point on the reduction if and only if, after using a matrix in 
$\GL_2(\OK)$ to replace $A$ and $B$ by suitable linear combinations,
we have $\pi^2 \dv a_{11}$, $\pi \dv a_{12}, a_{13}, a_{14}$ 
and $\pi \dv b_{11}$. 

\begin{Lemma}
Assume $\Char(k) \not= 2$ and let $Q_1 = Q_2 = 0$ be a minimal 
quadric intersection over $K$. We write $A$ and $B$ for the
$4$ by $4$ symmetric matrices corresponding 
to $Q_1$ and $Q_2$ and put $F(x,z) = \det (A x + B z)$. (If $Q_1=Q_2=0$
has level $0$ then the so-called doubling $y^2 = F(x,z)$ is again minimal.)
\begin{enumerate}
\item Suppose $(x:z) = (1:0)$ is a non-regular point on $y^2 = F(x,z)$
and let $s = 4 - \rank \widetilde{A}$. By a change of co-ordinates
we may assume
\begin{equation}
\label{eqn:AB}
  A = \begin{pmatrix} \pi A_1 & \pi A_2 \\ \pi A_2^T & A_3 
\end{pmatrix} \qquad B = \begin{pmatrix} B_1 & B_2 \\ B_2^T & B_3 
\end{pmatrix} 
\end{equation}
where $A_1$ and $B_1$ are $s$ by $s$ matrices. Let $q_1$ and $q_2$
be the quadratic forms corresponding to $A_1$ and $B_1$. Then
there are at most $4$ solutions to
\[ q_1(x_1,\ldots,x_s) \equiv q_2(x_1,\ldots,x_s) \equiv x_{s+1} \equiv 
\ldots \equiv x_4 \equiv 0
\pmod{\pi} \]
and each of these is a non-regular point on $Q_1= Q_2= 0$.
\item If we loop over all non-regular points on $y^2= F(x,z)$, moving each
to $(x:z) = (1:0)$ in turn, then all non-regular points on $Q_1= Q_2=0$ 
arise as described in (i).
\end{enumerate}
\end{Lemma}
\begin{Proof} (i) Since $Q_1 = Q_2 =0$ is minimal we have $s \le 3$.  
If $s=2$ then the binary quadratic forms
$q_1$ and $q_2$ cannot both vanish mod $\pi$ as this would contradict
minimality. Likewise if $s=3$ then $q_1$ and $q_2$ are ternary quadratic
forms with no common factor. So by  Bezout's theorem there are at most 
$4$ solutions. \\
(ii) Suppose $(1: 0: 0: 0)$ is a non-regular point. If we replace
$Q_1$ and $Q_2$ by suitable linear combinations then $A$ and $B$
are given by~(\ref{eqn:AB}) with $s= 1$ and $A_1 \equiv B_1 \equiv 
0 \pmod{\pi}$. 
It follows that $\det(Ax + Bz) 
= a x^4 + b x^3 z + \ldots$ with $\pi^2 \dv a$ and $\pi \dv b$. 
Then $(1:0)$ is a non-regular point on $y^2 = F(x,z)$. 
\end{Proof}

\begin{Remark}
These lemmas show that for a model of level $0$
the number of non-regular points is bounded independent of the size 
of the residue field.
This has the interpretation
that the $\OK$-scheme defined by the model is normal. 
Alternative proofs (taking a more geometric approach in the case $n=4$)
are given in~\cite{SadekThesis}. 
\end{Remark}

\subsection{Real solubility}
A section on testing local solubility would be incomplete without
some discussion of the real place. However we have nothing new to add.
For models of degree $3$ and for models of degree $2$ and $4$ 
with negative discriminant real solubility is automatic.
A binary quartic with positive discriminant has either $0$ or $4$
real roots, and in the former case is soluble over the reals if and
only if the leading coefficient is positive. For real solubility of
quadric intersections we refer to \cite[Chapter 6]{SiksekThesis}.
 
\section{Covering maps}
\label{sec:covmaps}

Let $\Phi$ be a non-singular genus one model over a field
$K$ with $\Char(K) \not= 2,3$. 
The starting point for this section is the survey article
\cite{AKM3P} that gives formulae for the covering map
$\pi : \CC_\Phi \to E$ where $E$ is the Jacobian elliptic curve
with Weierstrass equation~(\ref{eqn:jac}). 
The formulae are given by covariants coming from classical
invariant theory.

Our height bounds in Section~\ref{sec:htbds} will be computed
as sums of local contributions. To compute the correct 
contributions at primes dividing $2$ and $3$ we
modify the formulae in \cite{AKM3P}. The 
first step is to give a Weierstrass equation for the Jacobian 
\begin{equation}
\label{use-ainv}
 y^2 + a_1(\Phi) xy + a_3(\Phi) y = x^3 + a_2(\Phi) x^2 + a_4(\Phi) x 
+ a_6(\Phi) 
\end{equation}
that works in all characteristics. This is accomplished in \cite{ARVT},
\cite{CFS}, where the $a$-invariants $a_1, a_2, a_3, a_4, a_6$ 
are obtained from $c_4$ and $c_6$ by working back through the formulae
\begin{equation}
\label{formulaire}
\begin{aligned}
b_2 = a_1^2 +4 a_2, \qquad
b_4 = 2 &a_4 + a_1 a_3, \qquad
b_6 =  a_3^2 + 4 a_6, \\
c_4 = b_2^2-24 b_4,  \qquad
c_6 &=  -b_2^3 + 36 b_2 b_4 - 216 b_6.
\end{aligned}
\end{equation}
We recall formulae for the $a$-invariants below.
It is important to note however that they are
not invariants in the sense of Section~\ref{sec:g1m}.
Likewise our modified formulae for the covering maps will not 
be covariants. Nonetheless we still need to understand
how they change under transformations of genus one models.

\subsection{Generalised binary quartics}
\label{sec:covmap-bq}

We recall that a genus one model of degree~$2$ is a generalised
binary quartic $y^2 + P(x_1,x_2)y = Q(x_1,x_2)$ where
\begin{align*}
P(x_1,x_2) &= l x_1^2 + m x_1 x_2 + n x_2^2 \\
Q(x_1,x_2) &= a x_1^4 + b x_1^3 x_2 + c x_1^2 x_2^2 + d x_1 x_2^3 + e x_2^4.
\end{align*}
Let $g = \tfrac{1}{4}P^2 + Q$ be the binary quartic obtained by 
completing the square. It has covariants $h = \tfrac{1}{3} 
(g_{12}^2 - g_{11} g_{22})$ and $k= \tfrac{1}{12}(g_2 h_1 - g_1 h_2)$ 
where the subscripts denote partial derivatives.
In \cite{CFS} the $a$-invariants of $(P,Q)$ are defined as
\begin{align*}
a_1 &= m \\
a_2 &= -l n + c \\ 
a_3 &= l d + n b \\
a_4 &= -l^2 e - l n c - n^2 a - 4 a e + b d \\
a_6 &= -l^2 c e + l m b e - l n b d - m^2 a e 
        + m n a d - n^2 a c - 4 a c e + a d^2 + b^2 e.
\end{align*}
The $b$-invariants $b_2,b_4,b_6$ and $c$-invariants $c_4,c_6$ are 
then given by~(\ref{formulaire}).
We put $F = 4g = P^2 + 4Q$ and 
\begin{align*}
Z &= 2y + P \\
X &= \tfrac{1}{3}(h - b_2 g) \\
Y &= k - \tfrac{1}{2} a_1 XZ - \tfrac{1}{2} a_3 F Z.
\end{align*}

\begin{Lemma} (i) $Z, X, Y$ have coefficients in 
$\Z[l,m,n,a,b,c,d,e]$. \\
(ii) Let $(P,Q)$ be a non-singular 
generalised binary quartic defined over $K$.
Then $E = \Jac \CC_{(P,Q)}$ has Weierstrass equation
\[ y^2 + a_1 xy + a_3 y = x^3 + a_2 x^2 + a_4 x + a_6 \]
and the $2$-covering map $\CC_{(P,Q)} \to E$ is given by
$(x_1:x_2:y) \mapsto (X/Z^2,Y/Z^3)$.
\end{Lemma}
\begin{Proof}
(i) A direct calculation. \\
(ii) The formula for $E$ is recalled from \cite{CFS}.
The classical syzygy
\[ 27 k^2 = h^3 - 3 c_4 g^2 h - 2 c_6 g^3 \]
becomes
\begin{align*}
Y^2 + a_1 X Y Z + a_3  Y ZF &= X^3 + a_2 X^2 F + a_4 X F^2 + a_6 F^3 \\ 
& \qquad \qquad - (a_1 X + a_3 F)^2 (y^2 + P y - Q).
\end{align*} 
Since $F \equiv Z^2 \mod{(y^2 + Py - Q)}$ this gives the required map.
\end{Proof}

For use in later sections we put $F_2= F = P^2 + 4Q$ and $G_2 = X$.
Explicitly
\begin{align*}
  F_2 &= 
    (l^2 + 4 a) x_1^4 + (2 l m + 4 b) x_1^3 x_2 
    + (2 l n + m^2 + 4 c) x_1^2 x_2^2 
  + (2 m n + 4 d) x_1 x_2^3\\ &   + (n^2 + 4 e) x_2^4, \\
  G_2 &=
(-l^2 c + l m b - m^2 a - 4 a c + b^2) x_1^4 
+ (-2 l^2 d + 2 l n b - 4 m n a - 8 a d)
 x_1^3 x_2 \\ 
& + (-4 l^2 e - l m d + 2 l n c 
   - m n b - 4 n^2 a - 16 a e - 2 b d) x_1^2 x_2^2 \\
& + (-4 l m e + 2 l n d 
    - 2 n^2 b - 8 b e) x_1 x_2^3 + (-m^2 e +
    m n d - n^2 c - 4 c e + d^2) x_2^4.
\end{align*}
In \cite{SillsThesis} these polynomials were denoted $4G$ and $\widetilde{G}$.
We describe how they change under transformations
of genus one models. 
\begin{Lemma} 
\label{lem:FG2}
(i) If $(P',Q') = [\mu,(r_0,r_1,r_2),I_2](P,Q)$ then
\begin{align*}
F'_2(x,z) & = \mu^2 F_2(x,z) \\
G'_2(x,z) & = \mu^4 (G_2(x,z) 
   + (l r_2 + 2 r_0 r_2 + n r_0)F_2(x,z)). 
\end{align*}
(ii) If $(P',Q') = [1,0,(\begin{smallmatrix} \alpha & \beta \\ \gamma & \delta 
\end{smallmatrix}) ](P,Q)$ then
\begin{align*}
F'_2(x,z) & = F_2(\alpha x + \gamma z,\beta x + \delta z) \\
G'_2(x,z) & = (\alpha \delta - \beta \gamma)^2 G_2(\alpha x + \gamma z,\beta x + \delta z) 
  - \lambda F_2(\alpha x + \gamma z,\beta x + \delta z)
\end{align*}
where 
$ \lambda =  2 \alpha^2 \gamma^2 a 
  + \alpha \gamma (\alpha \delta+\beta \gamma) b
  + 2 \alpha \beta \gamma \delta c 
  + \beta \delta (\alpha \delta + \beta \gamma) d 
  + 2 \beta^2 \delta^2 e $. 
\end{Lemma}
\begin{Proof}
A direct calculation.
\end{Proof}

\subsection{Ternary cubics} 
A genus one model of degree 3 is a 
ternary cubic 
\begin{align*}
U(x_1,x_2,x_3) = a x_1^3 &+ b x_2^3 + c x_3^3 
   + f x_2^2 x_3 + g x_3^2 x_1 +  h x_1^2 x_2  \\
&  ~\text{~\qquad\qquad\qquad~}
  ~+ i x_2 x_3^2 + j x_3 x_1^2 + k x_1 x_2^2 + m x_1 x_2 x_3. 
\end{align*}
It has Hessian $H = -(1/2) \det (U_{ij})$ and covariants
\[  \Theta = (1/3) \left| \begin{array}{cccc} 
U_{11} & U_{12} & U_{13} & H_1 \\
U_{21} & U_{22} & U_{23} & H_2 \\
U_{31} & U_{32} & U_{33} & H_3 \\
H_{1} & H_{2} & H_{3} & 0 \end{array} \right|, 
\qquad
 J = (1/18) \left| \begin{array}{cccc} 
U_{1} & U_{2} & U_{3} \\
H_{1} & H_{2} & H_{3} \\
\Theta_{1} & \Theta_{2} & \Theta_{3} \end{array} \right|, 
\]
where the subscripts denote partial derivatives.
In \cite{ARVT}, \cite{CFS} the $a$-invariants of $U$ are defined as
\begin{equation*}
  \begin{aligned} \smallskip
    a_1 &=     m \\
  a_2 & =   -(f j + g k + h i) \\
  a_3 & =   9 a b c - a f i - b g j - c h k - f g h - i j k \\
  a_4 & =    -3 (a b g i + a c f k + b c h j) + a f^2 g + a i^2 k + b g^2 h + b i j^2 \\
& ~\qquad \qquad + c f h^2 + c j k^2 + f g j k + f h i j + g h i k \\
  a_6 & =    -27 a^2 b^2 c^2 + 9 a b c (a f i + b g j +c h k) + \ldots 
+ a b c m^3.
  \end{aligned}
\end{equation*}
The $b$-invariants $b_2,b_4,b_6$ and $c$-invariants $c_4,c_6$ are 
then given by~(\ref{formulaire}). We put
$b_8 = (b_2 b_6 - b_4^2)/4$ and
\begin{align*}
  Z &= \tfrac{1}{4} (H + b_2 U) \\
  X &= \tfrac{1}{192} (\Theta - 16 b_2 Z^2 - 12 b_2^2 Z U + b_2^3 U^2) \\
  Y &= \tfrac{1}{2} (\tfrac{1}{384} J
         - (a_1 X Z + a_3 Z^3 + a_3 X U + a_1 b_6 Z U^2 + a_1 b_8 U^3)).
\end{align*}
\begin{Lemma} (i) $Z,X,Y$ have coefficients in 
$\Z[a,b,c,f,g,h,i,j,k,m]$. \\
(ii) Let $U$ be a non-singular ternary cubic defined over $K$.
Then $E = \Jac \CC_{U}$ has Weierstrass equation
\[ y^2 + a_1 xy + a_3 y = x^3 + a_2 x^2 + a_4 x + a_6 \]
and the $3$-covering map $\CC_{U} \to E$ is given by
$ (x_1 : x_2 : x_3) \mapsto (X/Z^2,Y/Z^3)$.
\end{Lemma} 
\begin{Proof}
(i) A direct calculation. \\
(ii) The formula for $E$ is recalled from \cite{ARVT}, \cite{CFS}.
The classical syzygy
\begin{align*}
  12 J^2 &= \Theta^3 - 3 c_4 \Theta H^4  - 2 c_6 H^6 
    - 9 c_4 \Theta^2 H U + 12 c_6 \Theta H^3 U + 21 c_4^2 H^5 U \\
  & \qquad  + 6 c_6 \Theta^2 U^2 + 9 c_4^2 \Theta H^2 U^2 - 72 c_4 c_6 H^4 U^2 
    - 24 c_4 c_6 \Theta H U^3 \\ & \qquad  + (27 c_4^3 + 64 c_6^2) H^3 U^3 
     + 9 c_4^3 \Theta U^4 - 48 c_4^2 c_6 H^2 U^4 + 9 c_4^4 H U^5 
\end{align*}
becomes 
\begin{align*}
  Y^2 &+ a_1 X Y Z + a_3 Y Z^3 = X^3 + a_2 X^2 Z^2 + a_4 X Z^4 + a_6 Z^6 \\
  & - a_3 X Y U + (4 a_1 a_3 + 9 a_4) X^2 Z U 
  + \gamma_1 X Z^3 U
+ \gamma_2 Z^5 U 
   - (7 a_3^2 + 27 a_6) X^2 U^2 \\ & - (a_1 a_3^2 + 4 a_1 a_6) Y Z U^2 
+ \gamma_3 X Z^2 U^2 
+ \gamma_4 Z^4 U^2
+ \gamma_5 Y U^3 
+\gamma_6 X Z U^3 \\ &
+\gamma_7 Z^3 U^3 
+ \gamma_8 X U^4 
+ \gamma_9 Z^2 U^4 
+ \gamma_{10} Z U^5
+ \gamma_{11} U^6
\end{align*} 
where the $\gamma_i$ are certain polynomials in $\Z[a_1,a_2,a_3,a_4,a_6]$.
This gives the required formula for the $3$-covering map.
\end{Proof}

For use in later sections we put $F_3= Z^2$ and $G_3 = X$.
We describe how these polynomials change under transformations
of genus one models. 
\begin{Lemma} 
\label{lem:FG3}
(i) If $U' = [\mu,I_3]U$ then $F'_3 = \mu^6 F_3$ and $G'_3 = \mu^8 G_3$. \\
(ii) If $U' = [1,N ]U$ and $x_j = \sum n_{ij} x'_i$ where $N=(n_{ij})$ then 
\begin{equation}
\label{FG3}
\begin{aligned}
F'_3(x'_1,x'_2,x'_3) & = ((\det N)^4 F_3 + \alpha Z U + 
  \beta U^2)(x_1,x_2,x_3) \\
G'_3(x'_1,x'_2,x'_3) & = ((\det N)^6 G_3 + \lambda F_3 + 
  \gamma Z U + \delta U^2)(x_1,x_2,x_3) 
\end{aligned}
\end{equation}
for some $\lambda,\alpha,\beta,\gamma,\delta 
\in \Z[n_{11},n_{12}, \ldots,n_{33}, a,b,c, \ldots ,m]$. Moreover
if $N$ is diagonal then $\lambda=0$. 
\end{Lemma}
\begin{Proof}
(i) This is clear. \\
(ii) Since $H$ and $\Theta$ are covariants we have
\begin{align*} 
H'(x'_1,x'_2,x'_3) & = (\det N)^2 H(x_1,x_2,x_3) \\
\Theta'(x'_1,x'_2,x'_3) & = (\det N)^6 \Theta(x_1,x_2,x_3). 
\end{align*}
Let $\xi = \tfrac{1}{12} (b'_2 - (\det N)^2 b_2)$. Then~(\ref{FG3}) 
holds with $\lambda = - (\det N)^4 \xi$ and
\begin{align*} 
\alpha &= 6 (\det N)^2 \xi, & 
\beta &= 9 \xi^2, \\
\gamma &= -(\det N)^2 \xi (2 b'_2 - 9 \xi), &
\delta &=  -3 \xi^2 (b'_2 - 3 \xi).
\end{align*}
A generic calculation shows that $b'_2 \equiv (\det N)^2 b_2 \pmod{12}$.
Moreover if $N$ is diagonal then $b'_2 = (\det N)^2 b_2$ and so in that case
$\lambda = 0$.
\end{Proof}

\subsection{Quadric intersections}
\label{sec:qi}

A genus one model of degree $4$ is a pair of quadratic forms
$(Q_1,Q_2)$ in variables $x_1, \ldots, x_4$. We write
\begin{align*}
Q_1(x_1, \ldots, x_4) &= \sum_{i \le j} a_{ij} x_i x_j 
= \tfrac{1}{2} \sum_{i,j=1}^4 A_{ij} x_i x_j \\
Q_2(x_1, \ldots, x_4) &= \sum_{i \le j} b_{ij} x_i x_j 
= \tfrac{1}{2} \sum_{i,j=1}^4 B_{ij} x_i x_j 
\end{align*}
where $A = (A_{ij})$ and $B = (B_{ij})$ are the matrices of second
partial derivatives of $Q_1$ and $Q_2$. Let $Q^*_1 = \sum_{i \le j} 
a^*_{ij} x_i x_j$ and $Q^*_2 = \sum_{i \le j} b^*_{ij} x_i x_j$
be the quadrics whose matrices of second partial derivatives are
$\adj A$ and $\adj B$. There are covariants
\begin{align*}
T_1 &= \sum_{i,j=1}^4  \sum_{r \le s} b_{rs}^* (A_{ij}
A_{rs} - A_{is} A_{jr} ) x_i x_j \\
T_2 &= \sum_{i,j=1}^4  \sum_{r \le s} a_{rs}^* (B_{ij}
B_{rs} - B_{is} B_{jr} ) x_i x_j \\
J & =  (1/4) \, \frac{\partial (Q_1,Q_2,T_1,T_2)}
{\partial(x_1,x_2,x_3,x_4)}. 
\end{align*}

It is noted in \cite{CFS} that if 
$\Gamma = \sum_{i \le j} c_{ij} x_i x_j$ is a quadric in  4 variables then  
 \[ \det \left( \frac{\partial^2 \Gamma}{\partial x_i \partial x_j} \right)
  = \pf(\Gamma)^2 + 4 \rdet(\Gamma) \]
where $\pf(\Gamma) = c_{12} c_{34} + c_{13} c_{24} + c_{14} c_{23}$ and 
$\rdet(\Gamma) \in \Z[c_{11},c_{12}, \ldots, c_{44}]$.  
Writing $\pf (x Q_1 + z Q_2) = l x^2 + m x z + n z^2$ we put
\begin{align*} 
  Y &= \tfrac{1}{2} \left( J - l T_1^2 + m T_1 T_2 - n T_2^2 
         + m n (l T_1 + m T_2) Q_1  + l m (m T_1 + n T_2) Q_2 \right. \\
&~ \left. \hspace{10em}  
  + l^2 n^3 Q_1^2 +  l m n (l n + m^2) Q_1 Q_2 + l^3 n^2 Q_2^2 
\right).
\end{align*}
\begin{Lemma}
\label{pi4}
(i) $T_1,T_2,Y$ have coefficients in 
$\Z[a_{11}, a_{12}, \ldots , b_{44}]$. \\
(ii) Let $(Q_1,Q_2)$ be a non-singular quadric intersection 
defined over $K$. Then $(P,Q) = 
(\pf (x Q_1 + z Q_2),\rd(x Q_1 + z Q_2))$ is a non-singular 
generalised binary quartic 
and the $4$-covering map $\CC_{(Q_1,Q_2)} \to E = \Jac \CC_{(Q_1,Q_2)} $ 
is the composite of 
\[  \CC_{(Q_1,Q_2)} \to \CC_{(P,Q)} \,; \quad (x_1: \ldots : x_4) 
\mapsto (T_1:-T_2:Y) \]
and the $2$-covering map  $\CC_{(P,Q)} \to E$.
\end{Lemma} 
\begin{Proof}
(i) A direct calculation. \\
(ii) The formula for $(P,Q)$ is recalled from \cite{CFS}. There
is a classical syzygy satisfied by $Q_1, Q_2, T_1, T_2, J$ 
and the coefficients of 
\begin{equation}
\label{quartic}
 F(x,z) = \det (A x + B z). 
\end{equation}
Setting $Q_1=Q_2=0$ it reduces to
$J^2 \equiv F(T_1,-T_2) \mod{(Q_1,Q_2)}$.
We have $F = P^2 + 4Q$ and $2Y \equiv J - P(T_1,-T_2) \mod{(Q_1,Q_2)}$.
Therefore
\begin{equation*}
4( Y^2 + P(T_1,-T_2) Y - Q(T_1,-T_2)) = S_1 Q_1 + S_2 Q_2 
\end{equation*}
for some $S_1, S_2$ in 
$\Z[a_{11}, a_{12}, \ldots , b_{44}][x_1, \ldots, x_4]$.
Since the generic quadrics $Q_1$ and $Q_2$ are coprime mod $2$
a similar identity holds without the factor of $4$. Hence 
\begin{equation*}
 Y^2 + P(T_1,-T_2) Y \equiv Q(T_1,-T_2) \mod{(Q_1,Q_2)}
\end{equation*}
as required.
\end{Proof}

The $a$-invariants of $(Q_1,Q_2)$ are defined to be
the $a$-invariants of $(P,Q)$. 
The transformations of genus one models defined 
in Section~\ref{sec:g1m} have 
the following effect on $(P,Q)$ and on $T_1$ and $T_2$.

\begin{Lemma}
\label{lem:double}
If $(Q'_1,Q'_2) = [M,N](Q_1,Q_2)$ then $(P',Q')=[\det N,r,M](P,Q)$
for some $r=(r_0,r_1,r_2)$ where the $r_i$ are integer coefficient
polynomials in the entries of $M$ and $N$ and the coefficients of 
$Q_1$ and $Q_2$. Moreover if $N$ is diagonal then $r=0$.
\end{Lemma}
\begin{Proof}
If 
$N = I_4$ then the result is clear. 
So suppose $(Q'_1,Q'_2) = [I_2,N](Q_1,Q_2)$. We must show that
\begin{align*}
P'(x,z) &= (\det N) P(x,z) + 2r(x,z) \\
Q'(x,z) &= (\det N)^2 Q(x,z) - (\det N) P(x,z) r(x,z) - r(x,z)^2 
\end{align*}
for some $r(x,z)= r_0 x^2 + r_1 xz + r_2 z^2$ where the $r_i$ are
integer coefficient polynomials in the entries of $N$ and the 
coefficients of $Q_1$ and $Q_2$.
But in characteristic~$2$ we recognise $P(x,z) = \pf(x Q_1 + z Q_2)$ as 
the Pfaffian of a skew-symmetric matrix. This gives the formula for $P'$.
The formula for $Q'$ follows since $P'^2 + 4Q'= (\det N)^2(P^2 + 4Q)$.
Moreover if $N$ is diagonal then $P'(x,z) = (\det N) P(x,z)$ and so in
that case $r=0$.
\end{Proof}

\begin{Lemma}
\label{lem:T1T2}
(i) If $(Q'_1,Q'_2) = [(\begin{smallmatrix} \alpha & \beta \\ \gamma & \delta 
\end{smallmatrix}),I_4](Q_1,Q_2)$ then
\begin{equation}
\label{T1T2}
\begin{aligned}
 T'_1 &= (\alpha \delta-\beta \gamma)^2(\delta T_1 + \gamma T_2)
+ \nu_1 Q_1+\nu_2 Q_2 \\
 T'_2 &= (\alpha \delta-\beta \gamma)^2(\beta T_1 + \alpha T_2) 
+ \nu_3 Q_1+\nu_4 Q_2
\end{aligned}
\end{equation}
where the $\nu_i$ are integer coefficient polynomials in
$\alpha,\beta,\gamma,\delta$ and the coefficients of $Q_1$ and $Q_2$. \\
(ii) If $(Q'_1,Q'_2) = [I_2,N](Q_1,Q_2)$ and
$x_j = \sum n_{ij} x'_i$ where $N=(n_{ij})$ then
\[ T'_i(x'_1, \ldots,x'_4) = (\det N)^2 T_i(x_1, \ldots, x_4)\] 
for $i=1,2$.
\end{Lemma}
\begin{Proof}
(i) Let $a,b,c,d,e$ be the coefficients of (\ref{quartic}) and $a',b',c',d',e'$
their analogues for $(Q'_1,Q'_2)$. Direct calculation shows that~(\ref{T1T2}) 
holds with 
\begin{align*}
  \nu_1 &= \tfrac{1}{6} (  \gamma c' + 3 \alpha d' 
   - (\alpha \delta-\beta \gamma)^2 (\gamma c + 3 \delta d)) \\
  \nu_2 &=  \tfrac{1}{6} ( \delta c' + 3 \beta  d'
   - (\alpha \delta-\beta \gamma)^2 (\delta c + 3 \gamma b)) \\
  \nu_3 &=  \tfrac{1}{6} ( \alpha c' + 3 \gamma b'
   - (\alpha \delta-\beta \gamma)^2 (\alpha c + 3 \beta d)) \\
  \nu_4 &=  \tfrac{1}{6} ( \beta  c' + 3 \delta b' 
   - (\alpha \delta-\beta \gamma)^2 (\beta c + 3 \alpha b)).
\end{align*}
Writing $a',b',c',d',e'$ as polynomials in 
$\alpha,\beta,\gamma,\delta,a,b,c,d,e$ we find that
$\nu_1, \nu_2, \nu_3, \nu_4$
belong to $\Z[\alpha,\beta,\gamma,\delta,a,b,c,d,e]$.
These formulae are related to the 
covariance of the Hessian as defined in~\cite{g1hess}. \\
(ii) Let $M_1$ and $M_2$ be the matrices of second partial derivatives
of $T_1$ and $T_2$. Direct calculation shows that 
\[   \adj (\adj(A) x + \adj(B) z ) =  a^2 A x^3 + a M_1 x^2 z + e M_2 x z^2  
+ e^2 B z^3. \]
The covariance of $T_1$ and $T_2$ then follows from properties of the
adjugate.
\end{Proof}

For use in later sections we put 
$F_4= F_2(T_1,-T_2)$ and $G_4=G_2(T_1,-T_2)$ where $F_2$ and $G_2$ are
the polynomials associated to the model $(P,Q)$ in Lemma~\ref{pi4}(ii).

\subsection{A geometric observation}
\label{sec:geomobs}
Let $\Phi$ be a genus one
model of degree $n \in \{2,3,4\}$ over a field $K$. 
Let $E$ be the (possibly singular) curve defined by the
Weierstrass equation with coefficients the $a$-invariants 
of $\Phi$. The formulae in the last three sections 
define a map $\pi : \CC_{\Phi} \to E$. If $\Phi$
is non-singular then $\CC_{\Phi}$ is a smooth curve of genus one, $E$
is the Jacobian elliptic curve and $\pi$ is the $n$-covering map. 
However to understand what happens at primes of bad reduction we 
are also interested in singular models.

The composite $\CC_{\Phi} \stackrel{\pi}{\to} E \stackrel{x}{\to} \PP^1$
is given by $(F_n:G_n)$ where $F_n$ and $G_n$ are 
the homogeneous polynomials of degree $2n$ associated to $\Phi$.

\begin{Theorem} 
\label{lem:geomobs}
Let $\Phi$ be a genus one
model of degree $n \in \{2,3,4\}$ over a field $K$. 
Let $P \in \CC_{\Phi}$ say
$P = (x_1: x_2:y)$ or $(x_1: \ldots :x_n)$. Then 
$F_n(x_1,\ldots,x_n) = G_n(x_1,\ldots,x_n) = 0$ 
if and only if $P$ is singular or lies on 
a component of $\CC_\Phi$ of degree at most $n-2$.
\end{Theorem}
\begin{Proof} We split into the cases $n=2,3,4$.

\paragraph{{\bf Case $\mathbf{n=2}$}} The generalised binary quartic
\[ y^2 + (l x^2 + m x z + n z^2) y = a x^4 + b x^3 z + c x^2
z^2 + d x z^3 + e z^4. \]
has associated polynomials
\begin{align*} 
F_2(x,z) & = (l^2 + 4 a) x^4 + (2 l m + 4 b) x^3 z 
  + (2 l n + m^2 + 4 c) x^2 z^2 + \ldots  \\
G_2(x,z) & = (-l^2 c + l m b - m^2 a - 4 a c + b^2) x^4 + \ldots
\end{align*}
By Lemma~\ref{lem:FG2} we may assume that $P$ is the point $(x:z:y) = (1:0:0)$
and so $a=0$. Then $F_2(1,0) = G_2(1,0)= 0$ if and only if 
$l = b= 0$. This is the condition for $P$ to be a singular point.

\paragraph{{\bf Case $\mathbf{n=3}$}} 
A genus one model of degree $3$ is a ternary cubic
\begin{align*}
U(x_1,x_2,x_3) = a x_1^3 &+ b x_2^3 + c x_3^3 
   + f x_2^2 x_3 + g x_3^2 x_1 +  h x_1^2 x_2  \\
&  ~\text{~\qquad\qquad\qquad~}
  ~+ i x_2 x_3^2 + j x_3 x_1^2 + k x_1 x_2^2 + m x_1 x_2 x_3. 
\end{align*}
By Lemma~\ref{lem:FG3} we may assume that $P$ is the point $(x_1:x_2:x_3) = (1:0:0)$
and $a=h=0$. We compute
\begin{align*}
 F_3(1,0,0) &= j^4 k^2 \\
 G_3(1,0,0) &= b^2 j^6 - b j^5 k m + f j^5 k^2. 
\end{align*}
Thus $F_3(1,0,0)=G_3(1,0,0)=0$ if and only if $j=0$ or $b = k=0$.
These are the conditions that $P$ is either a singular point or 
lies on a line.

\paragraph{{\bf Case $\mathbf{n=4}$}} 
By Lemmas~\ref{lem:FG2},~\ref{lem:double} 
and~\ref{lem:T1T2} we may assume
that $P$ is the point $(1:0:0:0)$ and $\Phi = (Q_1,Q_2)$ takes the form 
\begin{align*}
Q_1(x_1, \ldots,x_4) & = \la x_1 x_3 + q_1(x_2,x_3,x_4) \\
Q_2(x_1, \ldots,x_4) & = \mu x_1 x_4 + q_2(x_2,x_3,x_4). 
\end{align*}
We compute $T_1(1,0,0,0) = \la^2 \mu^2 b_{22}$ and $T_2(1,0,0,0) = \la^2 \mu^2 a_{22}$.
If $\lambda \mu = 0$ or $a_{22} = b_{22}=0$ then $F_4(P)=G_4(P)=0$ and 
$P$ is a either a singular point or lies on a line. 
Otherwise we may assume that $\la=\mu=b_{22}=1$ and $a_{22}=0$. Then $P$ maps
to the point $(x:z:y) = (1:0:0)$ on the generalised binary quartic
\begin{align*}
  y^2  + (a_{24}&  x^2 + (a_{23} + b_{24}) x z + b_{23} z^2) y = 
 -(a_{23} a_{24} + a_{44}) x^3 z + \\
 & -(a_{23} b_{24} + a_{24} b_{23} - 
   a_{34} + b_{44}) x^2 z^2 - (a_{33} + b_{23} b_{24} - b_{34}) x z^3 - b_{33} z^4. 
\end{align*}
Our proof in the case $n=2$ shows that $F_4(P)=G_4(P)=0$ if and only if
$a_{24} = a_{44} = 0$. This is the condition for some quadric in the pencil
spanned by $Q_1$ and $Q_2$ (in fact it can only be $Q_1$) 
to factor as a product of two linear forms.
It is therefore also the condition for $P$ to lie on a conic.
\end{Proof}

\begin{Remark}
We suspect that some analogue of Theorem~\ref{lem:geomobs} 
holds for $n$-coverings
more generally. However our method of proof, using invariant theory
and explicit formulae, is unlikely to generalise to larger $n$.
\end{Remark}

\section{Height bounds}
\label{sec:htbds}

Let $E$ be an elliptic curve over a number field $K$. 
An $n$-descent calculation on $E$ computes equations for
the everywhere locally soluble $n$-coverings $\pi : \CC \to E$.
It is expected that a point $P \in \CC(K)$ will have
smaller height than its image in $E(K)$, and that therefore
searching on the covering curves makes it easier to find
generators for $E(K)$. 
Of course such an expectation can only be realised if our equations 
for $\CC$ are given relative to some reasonably good 
choice of co-ordinates.
In \cite{CFS} it is explained (at least over $K= \Q$) how to 
make such choices of co-ordinates when $n=2,3$ or $4$. 
We determine explicit height bounds in these cases.

\subsection{Local height bounds}
\label{sec:locglob}

Let $\Phi$ be a non-singular genus one model of degree $n \in \{2,3,4\}$
over a number field $K$. Let $M_K$, respectively $M^0_K$, be the
set of places, respectively finite places, of $K$. We write $K_v$ for
the completion of $K$ at $v \in M_K$ and normalise the 
absolute values $|\cdot|_v$ on $K_v$ so that 
the product formula holds.
The height of a point $P=(x_1: \ldots:x_n) \in \PP^{n-1}(K)$ is
\[ h(P) = \log \prod_{v \in M_K} \max ( |x_1|_v, \ldots, |x_n|_v). \]

Let $F_n$ and $G_n$ be the polynomials associated to $\Phi$ as defined
in Section~\ref{sec:covmaps}. For $v \in M_K$ we define
\begin{equation*}
\delta_v (\Phi) = \sup_{P \in \CC_{\Phi}(K_v)} 
\frac{ \max (| F_n(\x) |_v, |G_n(\x) |_v) }
{ \,\,\, \max( |x_1|_v, \ldots, |x_n|_v )^{2n} }
\end{equation*}
\begin{equation*}
\eps_v (\Phi) = \inf_{P \in \CC_{\Phi}(K_v)} 
\frac{ \max (| F_n(\x) |_v, |G_n(\x) |_v) }
{ \,\,\, \max( |x_1|_v, \ldots, |x_n|_v )^{2n} }
\end{equation*}
where $P= (x_1:x_2:y)$ or $(x_1: \ldots:x_n)$. These definitions
are independent of the scaling of the $x_i$ since $F_n$ and $G_n$ are
homogeneous of degree $2n$.

\begin{Theorem}
Let $\Phi$ be a non-singular genus one model over $K$.
\begin{enumerate}
\item For any $v \in M_K$ 
we have $0< \eps_v(\Phi) \le \delta_v(\Phi) < \infty$. 
\item If $v \in M^0_K$ and $\Phi$ is $v$-integral then 
$0< \eps_v(\Phi) \le \delta_v(\Phi) \le 1$.
\item If $v \in M^0_K$ and $\Phi$ has good reduction mod $v$ then
$\eps_v(\Phi) = \delta_v(\Phi) = 1$.
\item Let $h$ and $h_E$ be the heights on $\CC_\Phi$ and 
$E = \Jac(\CC_{\Phi})$ relative to $\CC_\Phi \to \PP^{n-1}$ 
and the Weierstrass equation (\ref{use-ainv}).
Let $\pi : \CC_\Phi \to E$ be the covering map.
Then for $P \in \CC_{\Phi}(K)$ we have
\begin{equation}
\label{bound-locglob}
 \qquad -\sum_{v} \log \delta_v(\Phi) \le  2n h(P) - h_E( \pi P) 
 \le -\sum_{v} \log \eps_v(\Phi). 
\end{equation}
\end{enumerate}
\end{Theorem}
\begin{Proof} (i) We are assuming that $\Phi$ is non-singular.
So by Theorem~\ref{lem:geomobs} there does not exist $P \in \CC_{\Phi}(K_v)$ 
with $F_n(P)=G_n(P)=0$. Since $\CC_{\Phi}(K_v)$ is compact it follows 
that $0< \eps_v(\Phi) \le \delta_v(\Phi) < \infty$. \\
(ii) Let $\OO_v$ be the valuation ring of $K_v$. If $\Phi$ has coefficients 
in $\OO_v$ then so do $F_n$ and $G_n$. We scale the $x_i$ 
so that $\max ( |x_1|_v, \ldots, |x_n|_v) = 1$.
Then $|F_n(\x)|_v \le 1$ and $|G_n(\x)|_v \le 1$. 
Hence $\delta_v(\Phi) \le 1$. \\
(iii) Again we scale the $x_i$ so that $\max ( |x_1|_v, \ldots, |x_n|_v) = 1$.
Then by Theorem~\ref{lem:geomobs} applied to the reduction 
of $\Phi$ mod $v$ we have
$\max (| F_n(\x) |_v, |G_n(\x) |_v)=1$. 
Hence $\eps_v(\Phi) = \delta_v(\Phi) =1$. \\
(iv) If $P \in \CC_{\Phi}(K)$, say $P= (x_1:x_2:y)$ or $(x_1: \ldots :x_n)$, 
then 
\[ h(P) =  \log \prod_{v \in M_K} \max ( |x_1|_v, \ldots, |x_n|_v ) \]
and
\[ h_E(\pi P) =  \log \prod_{v \in M_K} \max ( |F_n(\x)|_v, |G_n(\x)|_v ). \]
Taking logs in the definitions of $\delta_v(\Phi)$ and $\eps_v(\Phi)$ 
and summing over $v \in M_K$ gives the result. 
Notice that by (i) we are taking logs of positive
numbers, and by (iii) the sums are finite.
\end{Proof}

If $v \in M_K^0$ with uniformiser $\pi_v$ then
\begin{equation}
\label{finitecontrib}
  \delta_v(\Phi) = | \pi_v |_v^{2 \min A_v(\Phi)} \quad \text{ and } 
\quad  \eps_v(\Phi) = | \pi_v |_v^{2 \max A_v(\Phi)} 
\end{equation}
where $A_v(\Phi)$ is the set of {\em Tamagawa distances} defined
and computed in the next two sections. 
An alternative description of the Tamagawa distances 
in Section~\ref{sec:bounds} explains the choice of name.
The computation of $\delta_v(\Phi)$ and $\eps_v(\Phi)$ for $v$ a real
place is the subject of Section~\ref{sec:infinite}.

\subsection{Computing the Tamagawa distances}
\label{sec:finite}
Let $K$ be a finite extension of $\Q_p$ with ring of integers
$\OK$, maximal ideal $\pi \OK$, residue field $k$ and normalised discrete 
valuation $v : K^\times \to \Z$. The corresponding 
absolute value is $|x| = c^{-v(x)}$  for some constant $c>1$. 
Reduction mod $\pi$ will be denoted $x \mapsto \widetilde{x}$. 

Let $\Phi$ a non-singular genus one model over $K$ 
of degree $n \in \{2,3,4\}$. 
Let $F_n$ and $G_n$ be the polynomials depending on $\Phi$ as defined
in Section~\ref{sec:covmaps}.

\begin{Definition}
The set of {\em Tamagawa distances} $A = A(\Phi)$ is defined by
\[ \left\{ \frac{ \max (| F_n(\x) |, |G_n(\x) |) }
{ \,\,\, \max( |x_1|, \ldots, |x_n| )^{2n} } : P \in \CC_{\Phi}(K) \right\}
= \left\{ |\pi|^{2 \alpha} : \alpha \in A(\Phi) \right\}. \]
where $P = (x_1:x_2:y)$ or $(x_1: \ldots: x_n)$. In particular
$\CC_{\Phi}(K) \not= \emptyset$ if and only if $A(\Phi) \not= \emptyset$.
\end{Definition}

\begin{Definition}
A transformation of genus one models $g \in \G_n(K)$ is {\em integral}, 
respectively {\em diagonal}, if it satisfies the following conditions.
\[ \begin{array}{cccc}
n & g & \text{ integral } & \text{ diagonal } \\ \hline
2 & [ \mu, r, N] & \mu \in \OO_K^\times, r \in \OO_K^3, N \in \GL_2(\OK) & 
{\text{ $r=0$ and $N$ diagonal }} \\
3 & [ \mu, N] & \mu \in \OO_K^\times, N \in \GL_3(\OK) & 
N {\text{ diagonal }} \\
4 & [ M, N] & M \in \GL_2(\OK), N \in \GL_4(\OK) & 
{\text{ $M$ and $N$ diagonal. }} \\
\end{array} \]
\end{Definition}

The first part of the following theorem shows that if $\Phi$ and $\Phi'$ 
are $\OK$-equivalent then they have the same set of 
Tamagawa distances. The second part describes the effect 
of a diagonal transformation that preserves the level. 

\begin{Theorem} 
\label{FG-trans}
Let $\Phi$ and $\Phi'$ be genus one models over $\OK$ 
with $\Phi' = g \Phi$ for some $g \in \G_n(K)$, say
$g = [\mu,r,N]$, $[\mu,N]$ or $[M,N]$. 
Let $P \in \CC_{\Phi}(K)$, say $P = (x_1:x_2:y)$ or $(x_1: \ldots : x_n)$,
and $P' \in \CC_{\Phi'}(K)$, say 
$P' = (x'_1:x'_2:y')$ or $(x'_1: \ldots : x'_n)$,
with $x_j = \sum n_{ij} x'_i$ where $N = (n_{ij})$. If either
(i) $g$ is integral or (ii) $\det g \in \OO^\times_K$ and $g$ is diagonal
then 
\[ \max (| F'_n(\x')|, |G'_n(\x')|) = | \det N |^{-2} 
  \max (| F_n(\x)|, |G_n(\x)|). \]
\end{Theorem}
\begin{Proof}
Let $(r_2,s_2) = (2,4)$, $(r_3,s_3) = (6,8)$, $(r_4,s_4) = (12,14)$.
By Lemmas~\ref{lem:FG2}, \ref{lem:FG3}, \ref{lem:double} 
and~\ref{lem:T1T2} we have
\[ \left( \begin{array}{c} F'_n(x'_1, \ldots, x'_n) \\
 G'_n(x'_1, \ldots, x'_n) \end{array} \right) = (\det N)^{-2} 
\begin{pmatrix} (\det g)^{r_n} & 0 \\ \lambda & (\det g)^{s_n} 
\end{pmatrix} \left( \begin{array}{c} F_n(x_1, \ldots, x_n) \\
 G_n(x_1, \ldots, x_n) \end{array} \right) \]
for some $\lambda \in K$. These lemmas also show that (i) if $g$
is integral then $\lambda \in \OK$ and (ii) if $g$ is diagonal 
then $\lambda = 0$.
Taking absolute values gives the result.
\end{Proof}

We use Theorems~\ref{lem:geomobs} and~\ref{FG-trans} to modify our
local solubility algorithms in Section~\ref{sec:locsol} to give algorithms 
for computing the set of Tamagawa distances. Our presentation differs
from these earlier algorithms in that we do not restrict attention 
to (points whose reduction lies on) an affine piece until after the 
first iteration. For models of degrees $3$ and $4$ we use the 
subalgorithms in Section~\ref{sec:linesandconics} to compute the 
contributions from lines and conics.
The proof that our algorithms terminate (for $\Phi$ non-singular) 
is given in Section~\ref{sec:bounds}.

\begin{Algorithm}
\label{tdist:n=2}
{\tt TamagawaDistances(P,Q,Affine) } \\
{\tt INPUT:} A generalised binary quartic 
$\Phi = (P,Q)$ over $\OK$ and a boolean {\tt Affine}. \\
{\tt OUTPUT:} A finite set of non-negative integers $A$ such that
\[  \left\{ \frac{ \max (| F_2(\x) | , |G_2(\x) | ) }
{ \,\,\, \max( |x_1| , |x_2| )^{4} } : (x_1: x_2:y) \in 
\CC_{\Phi}(K)^\dagger \right\} = \{ |\pi|^{2 \alpha} : \alpha \in A \}  \]
where $\CC_{\Phi}(K)^\dagger = \{ R \in \CC_{\Phi}(K) : 
\widetilde{R} \in \Gamma \}$ and $\Gamma$ is the curve over $k$ defined by
\begin{align*} ~\qquad
\{ y^2 + \widetilde{P}(x_1,x_2) y = \widetilde{Q}(x_1,x_2) \} 
&\subset \PP(1,1,2) && \text{ if {\tt Affine} $=$  {\tt FALSE}} \\ 
\{ y^2 + \widetilde{P}(x,1) y = \widetilde{Q}(x,1) \} 
&\subset \Aff^2 && \text{ if {\tt Affine} $=$ {\tt TRUE}.} 
\end{align*}
\begin{enumerate}
\item Set $A = \emptyset$.
\item 
If there are smooth $k$-points on $\Gamma$ then set $A = \{0\}$.
\item Find all non-regular $k$-points on $\Gamma$.
Use an $\OK$-transformation to move each such point to $(x_1:x_2:y) = (0:1:0)$.
Then compute \[ A_1 = {\tt TamagawaDistances(P1,Q1,TRUE)} \] where
$ P_1(x_1,x_2) = \pi^{-1} P(\pi x_1,x_2) $, 
$ Q_1(x_1,x_2) = \pi^{-2} Q(\pi x_1,x_2) $, 
and set $A = A \cup \{ \alpha+1 : \alpha \in A_1 \}$.
\item Return $A$.
\end{enumerate}
\end{Algorithm}

\begin{Algorithm}
\label{tdist:n=3}
{\tt TamagawaDistances(U,Affine) } \\
{\tt INPUT:} A ternary cubic $U \in \OK[x,y,z]$ 
and a boolean {\tt Affine}. \\
{\tt OUTPUT:} A finite set of non-negative integers $A$ such that
\[  \left\{ \frac{ \max (| F_3(\x) | , |G_3(\x) | ) }
{ \,\,\, \max( |x_1| ,|x_2|, |x_3| )^{6} } : (x_1: x_2:x_3) \in 
\CC_{U}(K)^\dagger \right\} = \{ |\pi|^{2 \alpha} : \alpha \in A \}  \]
where $\CC_{U}(K)^\dagger = \{ P  \in \CC_{U}(K) : \widetilde{P} \in \Gamma\}$
and $\Gamma$ is the curve over $k$ defined by
\begin{align*} ~\qquad
\{ \widetilde{U}(x,y,z) = 0 \} 
&\subset \PP^2 && \text{ if {\tt Affine} $=$ {\tt FALSE}} \\ 
\{ \widetilde{U}(x,y,1) = 0 \} 
&\subset \Aff^2 && \text{ if {\tt Affine} $=$ {\tt TRUE}.} 
\end{align*}
\begin{enumerate}
\item Set $A = \emptyset$.
\item 
If $\Gamma$ contains an absolutely irreducible component of degree 2 or 3 
then set $A = \{0\}$.
\item Find all $k$-rational lines that are components of $\Gamma$ 
of multiplicity one. 
Compute the contribution $\alpha$ of each such line using 
Proposition~\ref{subalg-line:n=3} and put $A = A \cup \{\alpha\}$.
\item Find all non-regular $k$-points on $\Gamma$.
Use a transformation in  $\GL_3(\OK)$ to move each such point 
to $(0:0:1)$.
Then compute \[ A_1 = {\tt TamagawaDistances(U1,TRUE)} \] where
$ U_1(x,y,z) = \pi^{-2} U(\pi x, \pi y, z) $ 
and set $A = A \cup \{ \alpha+2 : \alpha \in A_1 \}$.
\item Return $A$.
\end{enumerate}
\end{Algorithm}

\begin{Algorithm}
\label{tdist:n=4}
{\tt TamagawaDistances(Q1,Q2,Affine) } \\
{\tt INPUT:} A quadric intersection $\Phi = (Q_1,Q_2)$ over $\OK$
and a boolean {\tt Affine}. \\
{\tt OUTPUT:} A finite set of non-negative integers $A$ such that
\[  \left\{ \frac{ \max (| F_4(\x) | , |G_4(\x) | ) }
{ \,\,\, \max( |x_1| , \ldots , |x_4| )^{8} } : (x_1: \ldots :x_4) \in 
\CC_{\Phi}(K)^\dagger \right\} = \{ |\pi|^{2 \alpha} : \alpha \in A \}  \]
where $\CC_{\Phi}(K)^\dagger = \{ P \in \CC_{\Phi}(K) : \widetilde{P}
\in \Gamma \}$ and $\Gamma$ is the curve over $k$ defined by
\begin{align*} ~\qquad
\{ \widetilde{Q}_1(x_1,\ldots,x_4) = \widetilde{Q}_2(x_1,\ldots,x_4) = 0 \}
&\subset \PP^3 && \text{ if {\tt Affine} = {\tt FALSE}} \\ 
\{ \widetilde{Q}_1(x,y,z,1) = \widetilde{Q}_2(x,y,z,1) = 0 \}
&\subset \Aff^3 && \text{ if {\tt Affine} = {\tt TRUE}.} 
\end{align*}
\begin{enumerate}
\item Set $A = \emptyset$.
\item 
If $\Gamma$ contains an absolutely irreducible component of degree 3 or 4 
then set $A = \{0\}$.
\item Find all $k$-rational lines and conics that are components of 
$\Gamma$ of multiplicity one.
Compute the contribution $\alpha$ of each such component using
Propositions~\ref{subalg-conic:n=4} and~\ref{subalg-line:n=4}
and put $A = A \cup \{\alpha\}$. 
\item Find all non-regular $k$-points on $\Gamma$.
Use a transformation in  $\GL_4(\OK)$ to move each such point to 
$(0:0:0:1)$ and a transformation in $\GL_2(\OK)$ to arrange that
$\frac{\partial Q_1}{\partial x_j}(0,0,0,1) \equiv 0 \pmod{\pi}$ 
for $1 \le j \le 4$ and $Q_1(0,0,0,1) \equiv 0 \pmod{\pi^2}$.
Then compute \[ \qquad A_1 = {\tt TamagawaDistances(Q1',Q2',TRUE)} \] where
\begin{align*}
\qquad  
Q'_1(x_1, \ldots,x_4) &= \pi^{-2} Q_1(\pi x_1, \pi x_2, \pi x_3, x_4) \\ 
Q'_2(x_1, \ldots,x_4) &= \pi^{-1} Q_2(\pi x_1, \pi x_2, \pi x_3, x_4)
\end{align*}
and set $A = A \cup \{ \alpha+3 : \alpha \in A_1 \}$.
\item Return $A$.
\end{enumerate}
\end{Algorithm}

\subsection{Contributions from lines and conics}
\label{sec:linesandconics}

Let $\Phi$ be a non-singular genus one model over $\OK$. Suppose
that the reduction of $\CC_{\Phi}$ mod $\pi$ contains a $k$-rational 
curve $C$ as a component of multiplicity one. (The multiplicity one 
condition is equivalent to requiring that all but finitely 
many $\kbar$-points on $C$ are smooth points on the reduction.)
Theorem~\ref{lem:geomobs} shows that if $C$ has degree $n-1$ or $n$
then the points $P \in \CC_\Phi(K)$ whose reduction is a 
smooth point on $C$ contribute $\alpha = 0$ to the set of Tamagawa distances.
In this section we determine the contributions in the remaining 
cases, namely when $n=3$ and $C$ is a line, and when $n=4$ and 
$C$ is a conic or line.

\begin{Proposition}
\label{subalg-line:n=3}
Let $U \in \OK[x,y,z]$ be a non-singular ternary cubic whose 
reduction contains a $k$-rational line $L$ as a component of multiplicity one. 
Then there is an integer $\alpha$ such that
\[  \frac{ \max (| F_3(\x) | , |G_3(\x) | ) }
{ \,\,\, \max( |x| ,|y|, |z| )^{6} }  = 
  |\pi|^{2 \alpha}  \]
for all $(x:y:z) \in \CC_U(K)$ whose reduction is a smooth point
on $L$. Moreover if $L$ is the line $\{x=0\}$ then $\alpha$ 
may be computed as follows.
\begin{enumerate}
\item Set $\alpha = 0$.
\item Replace $U$ by $\pi^{-1} U(\pi x,y,z)$ and let $\alpha = \alpha+1$.
\item Write $U(x,y,z) = f_0 x^3 + f_1(y,z) x^2 + f_2(y,z) x + f_3(y,z)$
where the $f_i$ are binary forms of degree $i$.
If $\widetilde{f_2} \dv \widetilde{f_3}$ say 
\[ f_3(y,z) \equiv (a y + b z) f_2(y,z) \pmod{\pi} \]
for some $a,b \in \OK$ then substitute $x \leftarrow x 
- a y - b z$ and go to Step~(ii).
\item Return $\alpha$.
\end{enumerate}
\end{Proposition}
\begin{Proof}
Writing
\[U(x,y,z) = f_0 x^3 + f_1(y,z) x^2 + f_2(y,z) x + f_3(y,z)\]
we are given that $v(f_3) \ge 1$ and $v(f_2)=0$.
If $P=(u:v:w) \in \CC_U(K)$ reduces to a smooth point on $L$
then $u \equiv 0$ and $f_2(v,w) \not\equiv 0 \pmod{\pi}$.
In Step (ii) we replace $P$ by $(\pi^{-1}u:v:w)$. The increase
of $\alpha$ by $1$ is justified by Theorem~\ref{FG-trans}(ii) with
$[\mu,N]  = [\pi^{-1},\Diag(\pi,1,1)]$. After this transformation
we still have $f_2(v,w) \not\equiv 0 \pmod{\pi}$ but now
\[U(x,y,z) \equiv f_2(y,z) x + f_3(y,z) \pmod{\pi}.\]
Hence $P$ reduces to a smooth point on the rational curve parametrised by
\[ (s:t) \mapsto  ( - \widetilde{f}_3(s,t): s  \widetilde{f}_2(s,t) : 
t \widetilde{f}_2(s,t) ) \]
If $\widetilde{f_2} \dv \widetilde{f_3}$ then this is a line and the
substitution in Step~(iii) moves the line to $\{x = 0\}$. We then return 
to Step~(ii). Otherwise we have a curve of degree 2 or 3 
and by Theorem~\ref{lem:geomobs} 
there is no further contribution to the Tamagawa distance.

We show in the next section that the algorithm terminates.
\end{Proof}

\begin{Proposition}
\label{subalg-conic:n=4}
Let $\Phi  = (Q_1,Q_2)$ be a non-singular quadric intersection over
$\OK$ whose reduction contains a $k$-rational conic $C$ 
as a component of multiplicity one. 
Then there is an integer $\alpha$ such that
\[  \frac{ \max (| F_4(\x) | , |G_4(\x) | ) }
{ \,\,\, \max( |x_1| , \ldots , |x_4| )^{8} }  = 
  |\pi|^{2 \alpha}  \]
for all $(x_1:\ldots:x_4) \in \CC_{\Phi}(K)$ whose reduction 
is a smooth point on $C$.
Moreover if $C$ is contained in the plane $\{x_1=0\}$ then 
$\alpha$ may be computed as follows.
\begin{enumerate}
\item Set $\alpha = 0$.
\item Make a $\GL_2(\OK)$-transformation so 
that $\widetilde{Q}_1$ vanishes on $\{x_1=0\}$. 
Replace $(Q_1,Q_2)$ by 
$(\pi^{-1} Q_1(\pi x_1,x_2,x_3,x_4), Q_2(\pi x_1,x_2,x_3,x_4))$
and let $\alpha = \alpha+1$.
\item Write $Q_i(x_1, \ldots,x_4)= 
\lambda_i x_1^2 + \ell_i(x_2,x_3,x_4)x_1+ q_i(x_2,x_3,x_4)$ for $i=1,2$.
If $\widetilde{q_1}$ belongs to the ideal 
generated by $\widetilde{\ell_1}$ and $\widetilde{q_2}$ say
\[q_1  \equiv (a_2 x_2 + a_3 x_3 + a_4 x_4) \ell_1 + b q_2 \pmod{\pi}\]
for some $a_2,a_3,a_4,b \in \OK$ then substitute $x_1 \leftarrow 
x_1 - (a_2 x_2 + a_3 x_3 + a_4 x_4)$ and go to Step (ii).
\item Return $\alpha$.
\end{enumerate}
\end{Proposition}

\begin{Proof}
To simplify the notation in the proof we first make 
a substitution in $x_2,x_3,x_4$ so that the conic $C$ is 
parametrised by $(s:t) \mapsto (0 : s^2 : s t : t^2 )$.

We write $Q_i = \lambda_i x_1^2 + \ell_i(x_2,x_3,x_4)x_1 + q_i(x_2,x_3,x_4)$
for $i=1,2$. After the $\GL_2(\OK)$-transformation in Step (ii) 
we have $v(q_1) \ge 1$. We put 
\[g(s,t) = \widetilde{\ell}_1(s^2,st,t^2).\] 
By the Jacobian criterion 
$(0:s^2:st:t^2)$ is a smooth point on the reduction if and only if
$g(s,t) \not= 0$.
Our hypothesis that $C$ has multiplicity one is therefore equivalent 
to the statement that $g$ is not identically zero. 

Suppose $P = (u_1: \ldots:u_4)$ reduces to a smooth 
point on $C$. Then (assuming $u_1, \ldots, u_4$ belong
to $\OK$ but not all to $\pi \OK$) we have
$\ell_1(u_2,u_3,u_4) \not\equiv 0 \pmod{\pi}$. 
In Step (ii) we replace $P$ by $(\pi^{-1}u_1 : u_2:u_3:u_4)$. 
The increase of $\alpha$ by $1$
is justified by Theorem~\ref{FG-trans}(ii) with $[M,N] = [\Diag(\pi^{-1},1),
\Diag(\pi,1,1,1)]$. This transformation changes neither $\ell_1$
nor $q_2$ but we now have
\begin{align*}
Q_1(x_1, \ldots, x_4) & \equiv x_1 \ell_1 (x_2,x_3,x_4) + q_1(x_2,x_3,x_4) 
\pmod{\pi} \\
Q_2(x_1, \ldots, x_4) & \equiv \hspace{7.6em} q_2(x_2,x_3,x_4) \pmod{\pi}. 
\end{align*}
Hence $P$ reduces to a smooth point on the rational curve parametrised by
\[ (s:t) \mapsto (-f(s,t): g(s,t) s^2 : g(s,t) st : g(s,t) t^2) \]
where $f(s,t) = \widetilde{q}_1(s^2,st,t^2)$. Since $g$ is not
identically zero this is a curve of degree $2$, $3$ or $4$.
If it has degree $2$ then in Step (iii) we move it to lie in 
the plane $\{x_1= 0\}$ and return to Step (ii). Otherwise we have
a curve of degree $3$ or $4$ and by Theorem~\ref{lem:geomobs} 
there is no further contribution to the Tamagawa distance.

We show in the next section that the algorithm terminates.
\end{Proof}

\begin{Proposition}
\label{subalg-line:n=4}
Let $\Phi = (Q_1,Q_2)$ be a non-singular quadric intersection over $\OK$
whose reduction contains a $k$-rational line $L$ as a component of 
multiplicity one. Then there is an integer $\alpha$ such that
\[  \frac{ \max (| F_4(\x) | , |G_4(\x) | ) }
{ \,\,\, \max( |x_1| , \ldots , |x_4| )^{8} }  = 
  |\pi|^{2 \alpha}  \]
for all $(x_1:\ldots:x_4) \in \CC_{\Phi}(K)$ whose reduction 
is a smooth point on $L$. Moreover if $L$ is the line 
$\{ x_1 = x_2 = 0 \}$ then $\alpha$ may be computed as follows.
\begin{enumerate}
\item Set $\alpha = 0$.
\item Replace $Q_i$ by $\pi^{-1} Q_i(\pi x_1,\pi x_2,x_3,x_4)$
for $i=1,2$ and let $\alpha = \alpha+2$.
\item Write $Q_1 = \sum_{i \le j} a_{ij} x_i x_j$ and
$Q_2 = \sum_{i \le j} b_{ij} x_i x_j$, and put
\[ \qquad  C = \begin{pmatrix} a_{13} & a_{23} \\ b_{13} & b_{23} \end{pmatrix},
\qquad D = \begin{pmatrix} a_{14} & a_{24} \\ b_{14} & b_{24} \end{pmatrix}. \]
Then compute $g(s,t) = \det(s \widetilde{C} + t \widetilde{D})$ and
\[  \qquad \left( \begin{array}{c} f_1(s,t) \\ f_2(s,t) \end{array} \right) 
= \adj (s \widetilde{C} + t \widetilde{D})  
  \left( \begin{array}{c} \widetilde{Q}_1(0,0,s,t) \\ \widetilde{Q}_2(0,0,s,t) 
\end{array} \right).  \]
\item If $g$ divides both $f_1$ and $f_2$ say
\begin{align*}
\qquad f_1(s,t) &=  (\widetilde{\lambda}_1 s 
  + \widetilde{\mu}_1 t) g(s,t)  \\
f_2(s,t) &=  (\widetilde{\lambda}_2 s 
  + \widetilde{\mu}_2 t) g(s,t)  
\end{align*} for some $\lambda_1,\lambda_2, \mu_1,\mu_2 \in \OK$ then 
substitute $x_3 \leftarrow x_3 + \lambda_1 x_1 + \lambda_2 x_2$
and $x_4 \leftarrow x_4 + \mu_1 x_1 + \mu_2 x_2$ and go to Step (ii).
\item If $f_1$, $f_2$ and $g$
have a common linear factor then solve for a linear form 
$\ell \in \OK[x_1, \ldots, x_4]$ with
\[ \qquad \widetilde{\ell} ( -f_1(s,t), -f_2(s,t), g(s,t)s, g(s,t)t)= 0. \]
Make a $\GL_4(\OK)$-transformation so that $\ell = x_1$. Then 
run the algorithm of Proposition~\ref{subalg-conic:n=4} on $(Q_1,Q_2)$ 
and add the answer to $\alpha$.
\item Return $\alpha$.
\end{enumerate}
\end{Proposition}

\begin{Proof}
By the Jacobian criterion $(0:0:s:t)$ is a smooth point on the
reduction if and only if $g(s,t) \not=0$, where $g$ is as defined in 
Step (iii). Our hypothesis that $L$ has multiplicity one is therefore 
equivalent to the statement that $g$ is not identically zero. 

Suppose $P=(u_1: \ldots : u_4) \in \CC_\Phi(K)$ reduces to a smooth point 
on $L$. Then (assuming $u_1, \ldots, u_4$ belong
to $\OK$ but not all to $\pi \OK$) we have 
$g(\widetilde{u}_3,\widetilde{u}_4) \not= 0$. In Step (ii) 
we replace $P$ by $(\pi^{-1}u_1:\pi^{-1}u_2:u_3:u_4)$. The increase
in $\alpha$ by $2$ is justified by Theorem~\ref{FG-trans}(ii) with $[M,N] = 
[ \Diag(\pi^{-1},\pi^{-1}), \Diag(\pi,\pi,1,1)]$.  Solving for the
first two co-ordinates of $\widetilde{P}$ in terms of the
last two we find it is a smooth point on the rational curve 
parametrised by
\[ (s:t) \mapsto (-f_1(s,t):-f_2(s,t):g(s,t)s: g(s,t)t). \]
Since $g$ is not identically zero this is a curve
of degree $1$, $2$ or $3$. These cases are treated in Steps (iv),(v)
and (vi).

We show in the next section that the algorithm terminates.
\end{Proof}

\subsection{Bounds on the Tamagawa distances}
\label{sec:bounds}
We recall from Section~\ref{sec:g1m} that the discriminant is a certain
polynomial in the coefficients of a genus one model. 
In this section we bound the Tamagawa distances in terms 
of the valuation of the discriminant. In particular this proves
that our algorithms terminate. We then give an alternative 
description of the Tamagawa distances.

\begin{Lemma}
\label{lem:minors}
Let $D = (d_{ij})$ be the $2$ by $5$ matrix over $\Z[l,m,n,a,b,c,d,e]$
whose entries are the coefficients of $F_2$ and $G_2$ as defined
in Section~\ref{sec:covmap-bq}. Then 
\[ \Delta = - 27 m_{15}^2 + 4 m_{14} m_{25} - m_{13} m_{35} \]
where $m_{ij} = d_{1i} d_{2j} - d_{1j} d_{2i}$. 
\end{Lemma}
\begin{Proof}
A direct calculation.
\end{Proof}

Our algorithms for computing the Tamagawa distances
(see Sections~\ref{sec:finite} and~\ref{sec:linesandconics})
only make transformations that preserve the level. 

\begin{Definition}
\label{defn:types}
Let $g \in \G_n(K)$ be a transformation of genus one models of
degree $n \in \{2,3,4\}$, say $g = [\mu,r,N], [\mu,N]$ or $[M,N]$. 
Then $g$ is a transformation of {\em type $r$} with $0<r<n$ 
if $\det(g) \in \OO_K^\times$ and the Smith normal
form of $N$ is $\Diag(I_{n-r}, \pi I_r)$.
\end{Definition}

We establish the following bounds on the Tamagawa distances.

\begin{Theorem}
\label{thm:bound}
Let $\Phi$ be a genus one model over $\OK$ of 
degree $n \in \{ 2,3,4\}$. Then the set of Tamagawa distances 
$A(\Phi)$ is bounded by
\[ \max A (\Phi) \le \left\{ \begin{array}{ll}
\tfrac{1}{2} v(\Delta) & \text{ if } n= 2 \\
v(\Delta) & \text{ if } n= 3 \\
2 v(\Delta) & \text{ if } n= 4 
\end{array} \right. \]
where $\Delta = \Delta(\Phi)$. Moreover if $v(\Delta) = 1$ 
then $A(\Phi) = \{0\}$. 
\end{Theorem}

\begin{Proof} We split into the cases $n=2,3,4$.
\paragraph{{\bf Case $\mathbf{n=2}$}}
Let $y^2 + P(x,z) y = Q(x,z)$ be a generalised binary quartic with
coefficients $l,m,n$ and $a,b,c,d,e$.
By Lemma~\ref{lem:minors} the discriminant $\Delta$ belongs to the 
ideal $(n^2,nd,d^2,e)$ in $\Z[l,m,n,a,b,c,d,e]$. 
But if $\alpha$ is a Tamagawa distance then $(P,Q)$
is $\OK$-equivalent to a model with $\pi^\alpha \dv n,d$ and $\pi^{2 \alpha}
\dv e$. Hence $\pi^{2 \alpha} \dv \Delta$ and 
$\alpha \le  \frac{1}{2}v(\Delta)$.

\paragraph{{\bf Case $\mathbf{n=3}$}} 
We label the coefficients of our ternary cubic as
\begin{align*}
U(x_1,x_2,x_3) = a x_1^3 &+ b x_2^3 + c x_3^3 
   + f x_2^2 x_3 + g x_3^2 x_1 +  h x_1^2 x_2  \\
&  ~\text{~\qquad\qquad\qquad~}
  ~+ i x_2 x_3^2 + j x_3 x_1^2 + k x_1 x_2^2 + m x_1 x_2 x_3. 
\end{align*}
Let $I_1 = (a ,h,k,b)$ and $I_2 = (b,f,i,c)$ in 
$\Z[a,b,c, \ldots, m]$. We checked using Magma
that the discriminant $\Delta$ belongs to $I_1 I_2^2$.

Let $\alpha$ be a Tamagawa distance. 
Then $\alpha = \alpha_1 + 2 \alpha_2$ 
where Algorithm~\ref{tdist:n=3} performs $\alpha_r$ transformations of
type $r$. The ternary cubic passed to the subalgorithm in 
Proposition~\ref{subalg-line:n=3} is $\OK$-equivalent to one with
$\pi^{\alpha_1} \dv a ,h,k,b$ and $\pi^{\alpha_2} \dv b,f,i,c$.
Since $\Delta \in I_1 I_2^2$ it follows that
$\alpha = \alpha_1 + 2 \alpha_2 \le v(\Delta)$.
By symmetry we also have $\Delta \in I_1^2 I_2$ and so
$\alpha_1,\alpha_2 \le \frac{1}{2} v(\Delta)$. In particular 
if $v(\Delta) = 1$ then $\alpha=0$.

\paragraph{{\bf Case $\mathbf{n=4}$}} In Section~\ref{sec:qi} we saw that 
the quadric intersection $(Q_1,Q_2)$ 
has the same discriminant as the generalised binary quartic
\begin{equation}
\label{eqn:gbq}
 y^2 + \pf (x Q_1 + z Q_2)y  = \rd(x Q_1 + z Q_2). 
\end{equation}
As usual we label the coefficients $l,m,n$ and $a,b,c,d,e$. 
Applying Lemma~\ref{lem:minors} to this generalised binary quartic, 
the discriminant $\Delta$ belongs to 
$J_1J_2$ where $J_1 = (n^2,nd,d^2,e)$ and $J_2$ is the 
ideal generated by the $2 \times 2$ minors of $D$.

Let $\alpha$ be a Tamagawa distance. 
Then $\alpha = \alpha_1 + 2 \alpha_2 + 3 \alpha_3$ 
where Algorithm~\ref{tdist:n=4} performs $\alpha_3$ transformations of
type $3$, then $\alpha_2$ transformations of type $2$ and then
$\alpha_1$ transformations of type $1$. 
Notice that a transformation
of type $r$ has inverse of type $4-r$. 
The quadric intersection passed to the subalgorithm in 
Proposition~\ref{subalg-conic:n=4}
is both $\OK$-equivalent to a model $(Q_1,Q_2)$ with
\[ Q_2(0,x_2,x_3,x_4) \equiv 0 \pmod{\pi^{\alpha_1}}, \]
and $\OK$-equivalent to a model $(Q'_1,Q'_2)$ with
\[ Q'_1(x_1,x_2,0,0) \equiv Q'_2(x_1,x_2,0,0) \equiv 0 \pmod{\pi^{\alpha_2}}. \]
We may therefore assume that $\pi^{\alpha_1} \dv n,d$ and 
$\pi^{2 \alpha_1} \dv e$, and (using Lemma~\ref{lem:double} to check the 
conclusion is unaffected by an $\OK$-equivalence) 
that~(\ref{eqn:gbq}) is reducible 
mod $\pi^{\alpha_2}$, i.e. there are binary quadratic forms 
$t_1$ and $t_2$ satisfying
\begin{align*}
\pf (x Q_1 + z Q_2) & \equiv t_1(x,z) + t_2(x,z) \pmod{\pi^{\alpha_2}} \\
\rd (x Q_1 + z Q_2) & \equiv -t_1(x,z) t_2(x,z)  \pmod{\pi^{\alpha_2}}. 
\end{align*}
This last condition implies that 
the $2$ by $2$ minors of the matrix $D$
in Lemma~\ref{lem:minors} vanish mod $\pi^{\alpha_2}$. Since
$\Delta \in J_1J_2$ it follows that
$2 \alpha_1 + \alpha_2 \le v(\Delta)$. The same argument gives
$2 \alpha_3 + \alpha_2 \le v(\Delta)$.
Hence $\alpha = \frac{1}{2}(2 \alpha_1 + \alpha_2) 
+ \frac{3}{2}(2 \alpha_3 + \alpha_2) \le 2 v(\Delta)$.
By Lemma~\ref{lem:minors} we also have $\Delta \in J_2^2$
and so $\alpha_1,\alpha_2,\alpha_3 \le \frac{1}{2} v(\Delta)$.
In particular if $v(\Delta) = 1$ then $\alpha = 0$.

We have shown in the cases $n=2,3,4$ that if $v(\Delta) =1$ then 
$A(\Phi) \subset \{0\}$. To prove equality it remains to show that
any such model is $K$-soluble. Since $v(\Delta) =1$ we have
$v(\Delta_E)=1$ and so by Tate's algorithm the Tamagawa number
$c(E)$ is also $1$. By Lemma~\ref{lem:kerres} it suffices to prove 
$K^{\nr}$-solubility and this follows by the results in \cite{CFS}. 
\end{Proof}

\begin{Corollary}
When the input is a non-singular genus one model 
the algorithms in Sections~\ref{sec:algs}, 
\ref{sec:finite} and~\ref{sec:linesandconics} terminate.
\end{Corollary}
\begin{Proof}
For the algorithms in Sections~\ref{sec:finite} and~\ref{sec:linesandconics}
this is immediate from our bounds on the Tamagawa distances. 
Taking into account the transformations in Step (i) that 
immediately follow each recursion, the algorithms in Section~\ref{sec:algs} 
never increase the level. So after finitely many iterations 
the level is preserved.
(In practice we first run the algorithms in \cite{CFS}, and so the level 
is always preserved.) The proof of Theorem~\ref{thm:bound} shows 
that thereafter the number of iterations (all of type $n-1$) is 
bounded by $\tfrac{1}{2} v(\Delta)$.
\end{Proof}

\begin{Remark}
If we think of the
algorithms as performing a tree search, then Theorem~\ref{thm:bound}
bounds the depth of the search, and Section~\ref{sec:nonreg} 
(on  non-regular points) bounds the breadth of the search.
From both points of view it is clearly desirable that we first minimise
our model using the algorithms in \cite{CFS}.
\end{Remark}

For the rest of this section we assume that $\Phi$ is $K$-soluble 
and of level $0$. The set of Tamagawa distances $A(\Phi)$ 
has the following alternative interpretation.
Let $\NN$ be the set of all matrices $N$ in $\GL_n(K)$ such that for
some transformation $g = [\mu,r,N]$, $[\mu,N]$ or 
$[M,N]$ in $\G_n(K)$ the model $g \Phi$ is minimal (equivalently is 
integral of level $0$).
Let $\NN_0 \subset \NN$ be the subset where the reduction of $g \Phi$ 
defines a curve with a 
$k$-rational component of multiplicity one and degree $n-1$ or $n$. 
Let $G$ be the subgroup of $\GL_n(K)$ generated 
by $\GL_n(\OO_K)$ and the scalar matrices. Then
\[ A(\Phi) = \{ v( \det N_i) : 1 \le i \le m \} \]
where $N_1, \ldots, N_m$ are a set of representatives for $G \setminus \NN_0$
scaled so that each $N_i$ has entries in $\OK$ not all in $\pi \OK$.
 
Theorem~\ref{thm:bound} shows that the set 
$G \setminus \NN_0$ is finite. Alternatively this follows by work
of Sadek \cite{SadekThesis} who computes $\#(G \setminus \NN)$.
If $n > 2$ then the same methods show that 
$\#(G \setminus \NN_0)$ is the Tamagawa number $c(E)$
of $E = \Jac(\CC_\Phi)$. This is still true when $n=2$ if
we adopt the convention that models of degree $2$
whose reduction mod $\pi$ have two $k$-rational components 
are counted twice.

It is natural to consider the graph with vertex set $G \setminus \NN$
and (directed) edges corresponding to the transformations of types 
$1,2, \ldots, n-1$. We recall that $c(E)$ is the number of $k$-rational 
components of the special fibre of the N\'eron model. For each such 
component there is a preferred vertex where the component is seen 
as a curve of degree $n-1$ or $n$. These vertices make up the
set $G \setminus \NN_0$. We may interpret $A(\Phi)$ as the set
of distances (weighted by type) from the vertex corresponding to $\Phi$
to each of these special vertices. This explains why we call $A(\Phi)$
the set of Tamagawa distances.

These graphs are investigated further in \cite{SillsThesis} with particular
attention given to the case $n=4$ and $E$ with multiplicative reduction.
These investigations suggest that the bounds in Theorem~\ref{thm:bound} 
are best possible.

\subsection{Calculation at the infinite place}
\label{sec:infinite}

Since our examples in Section~\ref{sec:examples} are over
$K= \Q$ we will only consider real places. (If $n=2$ then the
complex places are already treated in \cite{CPS}.)

Let $\Phi$ be a non-singular genus one model over $\R$ of
degree $n \in \{2,3,4\}$. We assume $\CC_\Phi(\R) \not= \emptyset$.
Let $F_n$ and $G_n$ be the polynomials associated to $\Phi$ as defined
in Section~\ref{sec:covmaps} and let $r \in \R$. In this
section we compute
\begin{equation*}
\delta (\Phi,r) = \sup_{P \in \CC_{\Phi}(\R)} 
\frac{ \max (| F_n(\x) |, |r F_n(\x) + G_n(\x) |) }
{ \,\,\, \max( |x_1|, \ldots, |x_n| )^{2n} }
\end{equation*}
\begin{equation*}
\eps (\Phi,r) = \inf_{P \in \CC_{\Phi}(\R)} 
\frac{ \max (| F_n(\x) |, |r F_n(\x) + G_n(\x) |) }
{ \,\,\, \max( |x_1|, \ldots, |x_n| )^{2n} }
\end{equation*}
where $P= (x_1:x_2:y)$ or $(x_1: \ldots:x_n)$. These definitions
are slightly more general than those in Section~\ref{sec:locglob} 
as previously we took $r=0$.

\begin{Proposition}
\label{prop:real}
We can compute $\delta (\Phi,r)$, respectively $\eps (\Phi,r)$,
by taking the maximum, respectively minimum, over all points 
$P \in \CC_\Phi(\R)$ satisfying one of the following conditions:
\begin{enumerate}
\item $P = (x_1: \ldots : x_n)$ with $x_i = \pm x_j$ for some $i \not= j$,
\item $F_n(P) = \pm (r F_n(P) + G_n(P))$, 
\item $n=2$ and $F_2(P)=0$,
\item according as $n=2,3,4$, 
\[ \frac{\partial f}{\partial x_i}(P)  = 0, \qquad 
 \frac{\partial(U,f)}{\partial(x_i,x_j)}(P)  = 0, \qquad
\frac{\partial(Q_1,Q_2,f)}{\partial(x_i,x_j,x_k)}(P) = 0, \]
where $f = F_n$ or $r F_n + G_n$ and $i,j,k$ are distinct.
\end{enumerate}
\end{Proposition}

\begin{Proof}
Since $\CC_\Phi(\R)$ is non-empty we may identify it as the real
locus of an elliptic curve. In particular it is isomorphic as a smooth
real manifold to either one or two copies of the circle $\R/\Z$.
We are asked to find the maxima and minima of a continuous 
real-valued function on this manifold. In (i) and (ii) we consider
the points where this function is not differentiable, and in (iii) and 
(iv) we consider the points where its derivative vanishes. 
We recall by Theorem~\ref{lem:geomobs} that there are no points 
$P \in \CC_{\Phi}$ with $F_n(P)= G_n(P)=0$. Condition (iii)
is needed since after completing the square 
$\CC_\Phi$ has equation $y^2= F_2(x_1,x_2)$.
\end{Proof}

We check that the set of points $P$ in Proposition~\ref{prop:real} 
is finite. In case (i) it suffices to note (by Bezout's theorem) 
that $\CC_\Phi$ has finite intersection with any hyperplane. In
cases (ii) and (iii) we recall that $(F_n:G_n)$ defines a non-constant 
morphism $\CC_{\Phi} \to \PP^1$ and therefore has finite fibres. 
If there were infinitely many points $P$ satisfying one of
the conditions in case (iv) then (after 
permuting the co-ordinates if necessary) we would have
\[  \la F_n + \mu G_n \equiv x_1^{2n} \pmod{I} \]
for some $(\lambda : \mu) \in \PP^1(\R)$, where $I= 0$, $(U)$, $(Q_1,Q_2)$
according as $n=2,3,4$.
In particular the form
\[  \frac{\partial(F_2,G_2)}{\partial(x_1,x_2)} 
\quad \text{ or } \quad 
\frac{\partial(U,F_3,G_3)}{\partial(x_1,x_2,x_3)} 
\quad \text{ or } \quad 
\frac{\partial(Q_1,Q_2,F_4,G_4)}{\partial(x_1,x_2,x_3,x_4)} 
\]
would be divisible by $x_1^{2n-1}$. However the invariant theory
in Section~\ref{sec:covmaps} shows that these forms meet $\CC_\Phi$ 
in distinct points: namely 
$\pi^{-1}( E[2] \setminus \{0\})$ in the case $n=2$ 
and $\pi^{-1}(E[2])$ in the cases $n=3,4$. This is
the required contradiction.

Proposition~\ref{prop:real} allows us to compute
$\delta (\Phi,r)$ and $\eps (\Phi,r)$ numerically. 
The case $n=2$ is already covered in~\cite{Siksek}, 
\cite{CPS}. See \cite[Section 2.5]{SillsThesis} for a worked example.
In the cases $n=3,4$ we use the Gr\"obner basis machinery in Magma.
In Section~\ref{sec:examples} we consider models over $\Q$, 
so the Gr\"obner bases can be computed exactly.

\section{Examples}
\label{sec:examples}

\subsection{Explicit bounds}

Let $E/\Q$ be an elliptic curve with global
minimal Weierstrass equation~(\ref{minw})
and discriminant $\Delta_E$.
Let $\CC = \CC_\Phi$ be an $n$-covering of $E$, where $\Phi$
is a non-singular genus one model of degree $n \in \{2,3,4\}$. 
We assume that $\CC(\Q_p) \not=\emptyset$ and $\Phi$ is minimal 
at all primes $p$. 
Therefore by \cite[Theorem~3.4]{CFS} we have $\Delta(\Phi) = \Delta_E$.
In particular
$\CC$ and $E$ have the same primes of bad reduction.

In Sections~\ref{sec:finite} and~\ref{sec:linesandconics} 
we computed a finite set of integers $A_p = A_p(\Phi)$ at each 
bad prime $p$.
The Weierstrass equations~(\ref{minw}) and~(\ref{use-ainv}) 
are related by a substitution
\[  x \leftarrow x + r \qquad y \leftarrow y + s x + t \]
for some $r,s,t \in \Z$. In Section~\ref{sec:infinite} we computed
the real contributions $\delta_\infty(\Phi,r)$ and $\eps_\infty(\Phi,r)$.
The height bounds $B_1$ and $B_2$ in (\ref{eqn1}) are now given by
\begin{align*}
B_1 &= -(1/2n) \log \delta_\infty(\Phi,r) 
+ (1/n) \sum_{p \dv \Delta_E} \min A_{p}(\Phi) \log p \\
B_2 &= -(1/2n) \log \eps_\infty(\Phi,r) 
+ (1/n) \sum_{p \dv \Delta_E} \max A_{p}(\Phi) \log p
\end{align*}
This follows from~(\ref{bound-locglob}) and~(\ref{finitecontrib}),
except that in changing our choice of Weierstrass equation
(from that given by the $a$-invariants to a standard one)
we must replace $G_n$ by $r F_n + G_n$. This makes no change at the 
finite places since $r \in \Z$. 

By Theorem~\ref{thm:bound} we need only sum over 
primes $p$ with $p^2 \dv \Delta_E$.

\subsection{A first example}
Let $E$ be the elliptic curve
$y^2 + y = x^3 - 41079 x - 2440008$
labelled 120267g1 in \cite{Cremonabook}. The primes of bad reduction
are $p=  3, 7, 23, 83$ with Kodaira symbols 
${\rm I}_4^*, {\rm I}_4, {\rm I}_1, {\rm I}_3$
and Tamagawa numbers $4, 4, 1, 3$.
The group $E(\Q)$ is free of rank $2$ generated by 
$(-106, 850)$ and $(-157, 373)$.

Among the coverings of $E$
computed using $n$-descent for $n=2,3,4$ we choose the following for
illustration.
\begin{align*}
\CC_2 : \qquad  &  y^2 + z^2 y = -5 x^4 - 171 x^3 z 
 + 78 x^2 z^2 + 216 x z^3 - 106 z^4 \\   
\CC_3 : \qquad & 
12 x^2 y - 9 x^2 z + 9 x y^2 - 12 x y z + 7 y^3 + 10 y^2 z - 17 y z^2 - 6 z^3 = 0 \\
\CC_4 : \qquad & \left\{ \begin{array}{rcl} 
x_1 x_2 + x_1 x_3 + 3 x_1 x_4 + x_2 x_3 - 4 x_2 x_4 + x_3^2 + 6 x_3 x_4 + 2 x_4^2 & = & 0 \\
3 x_1 x_3 + 3 x_1 x_4 - x_2^2 + x_2 x_3 - 9 x_3^2 + 4 x_3 x_4 + x_4^2 & = & 0 
\end{array} \right.
\end{align*}
The sets of Tamagawa distances $A_p$ are as follows.
We compute these as multisets so that, as a check on our calculations, 
the size of $A_p$ is equal to the Tamagawa number.
(See the comments at the end of Section~\ref{sec:bounds}.)
\begin{align*}
& \quad n = 2 && \quad n= 3 && \quad n=4 \\
A_{3}  &= \{ 0, 0, 1, 1 \} & A_{3}  &= \{ 2, 3, 3, 4 \} 
& A_{3}  &= \{ 2, 4, 6, 8 \} \\
A_{7}  &= \{ 0, 0, 1, 1 \} & A_{7}  &= \{ 1, 1, 1, 2 \} 
& A_{7}  &= \{ 1, 2, 3, 4 \} \\ 
A_{23}  &= \{ 0 \} & A_{23}  &= \{ 0 \} & A_{23}  &= \{ 0 \}  \\
A_{83}  &= \{ 0, 0, 1 \} & A_{83}  &= \{ 0, 1, 2 \} & A_{83}  
&= \{ 1, 2, 2 \} 
\end{align*}

Combining these with the contributions at the infinite place we obtain the
following bounds on the height of $P_n \in \CC_n(\Q)$ 
mapping down to $P \in E(\Q)$.
\begin{align*}
 -3.06805 \le & h(P_2) - \tfrac{1}{4} h_E(P) \le 1.21943 \\
-2.80610  \le & h(P_3) - \tfrac{1}{6} h_E(P) \le 2.44241  \\
-3.08885 \le & h(P_4) - \tfrac{1}{8} h_E(P) \le 2.48228 
\end{align*}

The curves $\CC_n$ have many small rational points. 
We list a few of these together with their contributions
to the Tamagawa distances (at $p=3,7,83$) 
and the height difference $h(P) - \tfrac{1}{2n} h_E(\pi P)$.
\[ \begin{array}{c@{\qquad}lccccrc}
& \multicolumn{1}{c}{P} & p=3 & p=7 & p=83 & 
\multicolumn{3}{c}{ h(P) - \tfrac{1}{2n} h_E(\pi P)} \\ \hline
& (1:1:3) & 0 & 1 & 0 && -1.68305 \\ 
& (2:3:37) & 0 & 0 & 1 && -0.30284  \\ 
n=2 & (6:-1:178) & 1 & 1 & 0 && -1.08967  \\ 
& (27:-1:871) & 1 & 1 & 1 && 1.14846  \\ 
& (769:787:2143781) & 0 & 0 & 0 && -2.63972 \\ \hline
& (1:0:0) & 4 & 1 & 0 && -1.15212 \\ 
& (1:-1:-1) & 3 & 1 & 0 && -2.16072 \\ 
n=3 & (2:-3:1) & 2 & 2 & 0 && -1.74660 \\ 
& (2:18:15) & 4 & 2 & 2 && 1.96488 \\ 
& (1:-6:20) & 2 & 1 & 0 && -2.38783  \\ \hline
& (1:0:0:0) & 4 & 2 & 1 && -0.70073 \\ 
& (-2:5:2:7) & 2 & 4 & 1 && -1.54491 \\ 
n=4 & (-3:3:1:8) & 6 & 3 & 1 && -0.80265  \\ 
& \text{\small $(557:544:-134:470)$} & 2 & 2 & 1 && -2.31493 \\ 
& \text{\small $(157397:2728:1502:-1438)$} & 8 & 3 & 2 &&  1.99552 \\ \hline
\end{array} \]

\medskip

\subsection{Searching for generators of large height}

We give two examples. The first is an example where the generator
was found by Michael Stoll using 4-descent  (see \cite[Section 7C]{CFS}).
The elliptic curve $E$ in the second example is taken from a list
of rank 1 curves (for which the generator had not been found) 
sent to us by Robert Miller. Although in both these examples
the elliptic curve has rank 1, the conductor is large
enough to make a Heegner point calculation difficult.

\begin{Example}
\label{ex1}
Let $E/\Q$ be the elliptic curve $y^2 = x^3 + 7823$. An $L$-value computation
shows that $\rank E(\Q) = 1$ and the generator is predicted to
have canonical height $h_1 = 77.61777 \ldots$ (if we assume 
$\Sha(E/\Q)$ is trivial).

Using the implementations of $2$-, $3$- and $4$-descent in Magma,
together with minimisation and reduction, we obtain the following
$n$-coverings of $E$.
\begin{align*}
\CC_2 : \qquad  & y^2 + (x^2 + z^2) y = -3 x^4 + 28 x^3 z
- 2 x^2 z^2 - 4 x z^3 + 10 z^4 \\   
\CC_3 : \qquad & x^3 + x^2 y - 4 x^2 z - 8 x y z + 8 x z^2 
   + y^3 - 5 y^2 z - 7 y z^2 + z^3 = 0 \\
\CC_4 : \qquad & \left\{ \begin{array}{rcl} 
2 x_1 x_2 + x_1 x_3 + x_1 x_4 + x_2 x_4 + x_3^2 - 2 x_4^2 & = & 0\\
x_1^2 + x_1 x_3 - x_1 x_4 + 2 x_2^2 - x_2 x_3 + 2 x_2 x_4 
- x_3^2 - x_3 x_4 + x_4^2 & = & 0
\end{array} \right.
\end{align*}
At each of the bad primes $p = 2,3,7823$ the elliptic curve $E$ 
has additive reduction with Kodaira symbol II. The finite primes 
make no contribution to our height bounds. If $P_n \in \CC_n(\Q)$ 
maps down to $P \in E(\Q)$ then our bounds work out as
\begin{align*}
-1.94921 \le & h(P_2) - \tfrac{1}{4} h_E(P) \le -0.92414 \\
-2.91485 \le & h(P_3) - \tfrac{1}{6} h_E(P) \le -1.41177 \\
-3.66288 \le & h(P_4) - \tfrac{1}{8} h_E(P) \le -2.43592
\end{align*}
The bounds established in~\cite{CPS} show that for $P \in E(\Q)$ we have
\[ -3.68143 \le h_E(P) - \widehat{h}_E(P) \le  0.74248 \]
where $\widehat{h}_E$ is the canonical height.
We write $P_n = (x_1: x_2:y)$, respectively $(x_1: \ldots : x_n)$,
where $x_1, \ldots, x_n$ are coprime integers. Taking 
$\widehat{h}_E(P) =  h_1$ we therefore expect to find 
$P_n \in \CC_n(\Q)$ with $H_n = \max(|x_1|, \ldots, |x_n|)$ 
in the following ranges. For comparison we list the actual 
points $P_n$.
\begin{align*}
15170781  & \le H_2 \le 
127792792 & 
P_2 &=( 10677130 : -42786483 : 5018494588774686 ) \\
12185  & \le H_3 \le 
114492 & 
P_3 &= ( 10445: -32922: 16423 ) \\
265 & \le H_4 \le 
1570 & 
P_4 &= ( 116: 207: 474: -332 )
\end{align*}
\end{Example}

Example~\ref{ex1} makes precise the statement 
that searching on an $n$-covering to find a generator for $E(\Q)$ 
becomes easier as $n$ increases.
For the actual searching we use the $p$-adic method due to 
Elkies~\cite{Elkies} and Heath-Brown, as implemented in Magma by 
Watkins. This takes time $O(H)$, respectively $O(H^{2/3})$, to search
for points of height up to $H$ on a $3$-covering, 
respectively $4$-covering.

\begin{Example}
\label{ex2}
Let $E_0$ be the elliptic curve $y^2 + x y + y = x^3 - x^2 - 2305 x + 43447$,
labelled 3850m1 in \cite{Cremonabook}, and $E$ the quadratic twist 
of $E_0$ by $d = -2351$. 
We fix a Weierstrass equation for $E$ of the form~(\ref{minw}).
The primes of bad reduction are
$p = 2, 5, 7, 11, 2351$ with
Kodaira symbols 
$ {\rm I}_1, {\rm II}^*, {\rm I}_2, {\rm I}_1, {\rm I}_0^* $
 and Tamagawa numbers
$ 1, 1, 2, 1, 2 $.
An $L$-value computation shows that
$\rank E(\Q) = 1$ and the generator is predicted to have canonical height
$h_1 = 182.01408 \ldots$ (if we assume $\Sha(E/\Q)$ is trivial). 
The torsion subgroup of $E(\Q)$ is trivial.

Using $4$-descent in Magma we obtain a $4$-covering $\CC_4$ 
of $E$ with equations
\small \begin{align*}
3 x_1^2 + 17 x_1 x_2 + x_1 x_3 + 7 x_1 x_4 - 5 x_2^2 + 11 x_2 x_3 
+ 6 x_2 x_4 + 5 x_3^2 + 9 x_4^2 &=  0 \\
10 x_1^2 + 7 x_1 x_2 - x_1 x_3 - x_1 x_4 + 4 x_2^2
 - x_2 x_3 - 13 x_2 x_4 + 14 x_3^2 - 30 x_3 x_4 + 18 x_4^2 & =  0 
\end{align*}  \normalsize
The Tamagawa distances for this quadric intersection 
are $A_{2}  = A_{11}  =  \{ 0 \}$, $A_{5}  = \{ 6 \}$, $A_{7}  = \{ 1, 1 \}$ 
and $A_{2351}  = \{ 4, 4 \}$. For $P_4 \in \CC_4(\Q)$ we obtain the bounds
\[  0.65550 \le h(P_4) - \tfrac{1}{8} h_E(\pi P_4) \le 0.94857. \]
The bounds in \cite{CPS} are now 
$-15.51194  \le h_E(P) - \widehat{h}_E(P) \le 8.73556 $.
We are therefore looking for $P_4 \in \CC_4(\Q)$ with
\[  21.46827 \le h(P_4) \le 24.79228. \]
A direct search is not practical.
We now explain how using the theory in this paper, we are 
nonetheless able to find this point. 
Our computation of the Tamagawa distances at $p=5$ and $p=2351$ 
suggests replacing $\CC_4$ by either $\CC'_4$
with equations
\small \begin{align*}
3 x_1^2 + 3 x_1 x_2 + 4 x_1 x_3 + 6 x_1 x_4 + 3 x_2^2 - 3 x_2 x_3 + 2 x_2 x_4 + 6 x_3^2 - 
28 x_3 x_4 + 11 x_4^2 & = 0 \\
4 x_1^2 + x_1 x_2 - 7 x_1 x_3 + 9 x_1 x_4 - 4 x_2^2 - 8 x_2 x_3 + 38 x_2 x_4 + 31 x_3^2 + 
14 x_3 x_4 + 16 x_4^2 & = 0 
\end{align*}  \normalsize
or $\CC''_4$ with equations
\small \begin{align*}
2 x_1^2 + 4 x_1 x_2 + 10 x_1 x_3 + 3 x_1 x_4 - 3 x_2^2 - 2 x_2 x_3 - 6 x_2 x_4 - 5 x_3^2 - 
10 x_3 x_4 - 21 x_4^2 & =  0 \\
14 x_1^2 + x_1 x_2 + 11 x_1 x_3 - 11 x_1 x_4 + 2 x_2^2 + 25 x_2 x_3 + 15 x_2 x_4 - 2 x_3^2 
- 24 x_3 x_4 + 12 x_4^2 & =  0. 
\end{align*}  \normalsize
Again we have reduced these models as described in \cite{CFS}.
We do not record the changes of co-ordinates used, since they may easily 
be recovered using the algorithm in~\cite{improve4}, as implemented in the 
Magma function {\tt IsEquivalent}. 

On $\CC'_4$ and $\CC''_4$ we have $A_{2} = A_{5} = A_{11} = \{0\}$,
$A_{7} = \{1,1\}$ and $A_{2351} = \{ 0, 4 \}$. So
the only finite primes to contribute to 
our height bounds are $p=7$ and $p=2351$. 
Moreover if we are willing to search on both curves then 
the contributions at $p=2351$ may be ignored.
Suppose $P_4 \in \CC_4(\Q)$, corresponds to 
$P'_4 \in \CC_4'(\Q)$ and $P''_4 \in \CC_4''(\Q)$, and
maps down to $P \in E(\Q)$. Then depending on the reductions of
these points mod $2351$, we have either
\begin{equation}
\label{bds1}
 -9.65955 \le h(P'_4) - \tfrac{1}{8} h_E(P) \le -9.29236 
\end{equation}
or
\begin{equation}
\label{bds2}
 -9.72818 \le h(P''_4) - \tfrac{1}{8} h_E(P) \le -9.35987. 
\end{equation}
Taking $\widehat{h}_E(P) = h_1$ it follows that either
\begin{equation}
\label{alt}
  11.15322  \le h(P'_4) \le 14.55134 \quad \text{ or } \quad  
  11.08459 \le h(P''_4) \le 14.48383. 
\end{equation}
If we are willing to search on only one of these curves 
then the upper bounds increase by $\log 2351 = 7.76259\ldots$. 

Magma's {\tt PointSearch} takes just a few seconds to 
find a point $P''_4 \in \CC_4''(\Q)$. We find 
the corresponding 
points $P_4 \in \CC_4(\Q)$ and $P'_4 \in \CC'_4(\Q)$ 
by making the relevant changes
of co-ordinates, and thus obtain
\begin{align*}
& P_4 = \text{\tiny $( -32083748086: 42638879317: 38411124781:
 22127244455 )$} & h(P_4) &= 24.47603 \ldots \\
& P'_4 = \text{\tiny $( 472320823 : 4111701909 : 
-2388802174 : -2139378517 )$}  & h(P'_4) &= 22.13710 \ldots \\
& P''_4 = ( 785047 : -840912 : 1542460 : -236990 ) &
h(P''_4) &= 14.24888 \ldots
\end{align*}
In particular we see it is
the second of the two possibilities in~(\ref{alt}) that holds.
These points map down to $P = (u/w^2,v/w^3) \in E(\Q)$ where
\begin{align*}
u & = 1757287936905025328253331560718272340242739349926447025094428588\backslash \\
  & \qquad 4833392724486595115 \\
v & = 4125077432494049001174441775597880344806917503465242447257595890\backslash \\
  & \qquad 83530835657373093470958302511042544245136026529511888663249   \\
w & = 364436547292608819468573335937957548482. 
\end{align*}
If $P_0 \in E(\Q)$ is a generator then (assuming we have carried out the 4-descent
rigorously) it lifts to a rational point on $\CC_4$. 
Combining our height bounds~(\ref{bds1}) and~(\ref{bds2}) 
with those in \cite{CPS} it follows that
\[  \widehat{h}_E(P_0) \ge 8 \times 9.29236 - 8.73556 = 65.60332. \]
Since $\widehat{h}_E(P) = 182.01408 \ldots$ we deduce (without the need for
any further searching) that $P$ is a generator for $E(\Q)$. 

\end{Example}

Example~\ref{ex2} shows the advantages of searching on several
different models of the same curve. One strategy would be
to search on $\prod_p c_p(E)$ models of each curve, so that
only the contributions to our height bounds at the infinite place 
are relevant. (These contributions do not appear to vary greatly between
the models, so long as we always reduce them.) However when 
$\prod_p c_p(E)$ is large then some compromise is needed and for
this the graphs in \cite{SillsThesis} are useful. Alternatively
it may be possible to adapt the $p$-adic point searching method
to search on several models of the same curve simultaneously.

\section*{Acknowledgements}
This article is based on the second author's PhD thesis.
We would like to thank Mohammad Sadek for many useful discussions.
All computer calculations in support of this work were performed using 
Magma \cite{Magma}.

\end{document}